%%%                    From Rolci Cipolatti
%%% 
%%%                    This a Plain TeX file
%%%
%%% *************************************************************
%%%                        BEGIN PREAMBLE
%%% *************************************************************
%\magnification=\magstep1
%\input pagina.tex
%\modelodepagina1
%\hoffset=1truecm
%\voffset=1truecm
\hsize=16truecm 
\vsize=21.6truecm
\baselineskip=12 true pt plus .5pt minus .5pt
\parskip=4pt
\parindent=10pt
\abovedisplayskip=14pt plus 3pt minus 9pt
\abovedisplayshortskip=0pt plus 3pt 
\belowdisplayskip=14pt plus 3pt minus 9pt
\belowdisplayshortskip=7pt plus 3pt minus 4pt
%\nopagenumbers
\font\Titulo=cmbx12
\font\titulosecao=cmbx10 scaled\magstephalf
\font\letrautor=cmcsc10
\font\letranota=cmr8
\font\letrasl=cmsl8
%%%
%%%  macros1
%%%
\font\TenEns=msbm10
\font\SevenEns=msbm7
\font\FiveEns=msbm5
\newfam\Ensfam
\def\Ens{\fam\Ensfam\TenEns}
\textfont\Ensfam=\TenEns
\scriptfont\Ensfam=\SevenEns
\scriptscriptfont\Ensfam=\FiveEns
\def\R{{\Ens R}}

\def\S{{\Ens S}}

\font\TenCM=cmr10
\font\SevenCM=cmr7
\font\FiveCM=cmr5
\newfam\CMfam
\def\CM{\fam\CMfam\TenCM}
\textfont\CMfam=\TenCM
\scriptfont\CMfam=\SevenCM
\scriptscriptfont\CMfam=\FiveCM

\def\ee{{\CM e\/}}
\def\unim{{\CM i}}

\setbox111=\hbox to 2truemm{\hrulefill}
\setbox222=\hbox to 2truemm{\vrule height 2truemm width .4truept\hfil
\vrule height 2truemm width .4truept}
\def\cqd{\vbox{\offinterlineskip\copy111\copy222\copy111}}
\def\mod#1{\vert #1\vert}
\def\norma#1#2{\Vert #1 \Vert_{#2}}
\def\Mod#1{\left\vert #1\right\vert}

\def\divergente{\mathop{\rm div}\nolimits}
\def\supp{\mathop{\rm supp}\nolimits}
\def\diam{\mathop{\rm diam}\nolimits}
\def\dist{\mathop{\rm dist}\nolimits}

\def\cc#1{\kern .7em \hfill #1\hfill\kern.7em}
\def\Lpp{L^p(Q)}\def\Hpp{\Lpp}

\def\Wpp#1{{\cal W}_{#1}}
\def\Wpptil#1{\widetilde{{\cal W}}_{#1}}

\def\cpl{\kern 1pt}
\def\build#1_#2^#3{\mathrel{\mathop{\kern 0pt#1}\limits_{#2}^{#3}}}

%%%
%%%  macros2
%%%
%  ---------------------------------
%  Macros for automatic numbering of formulas and theorems 
%  ---------------------------------

\newcount\numerosection
\newcount\eqnumer
\newcount\lemnumer
\newcount\fignumer
\numerosection=0
\eqnumer=0
\lemnumer=0
\fignumer=0

\def\numsection{\global\advance\numerosection by1
\global\eqnumer=0
\global\lemnumer=0
\global\fignumer=0
\the\numerosection}

\def\strutdepth{\dp\strutbox}
\def\marginalsigne#1{\strut
    \vadjust{\kern-\strutdepth\specialsigne{#1}}}
\def\specialsigne#1{\vtop to \strutdepth{
\baselineskip\strutdepth\vss\llap{#1 }\null}}

\font\margefont=cmr10 at 6pt
\newif\ifshowingMacros
\showingMacrosfalse

\def\cite#1{\csname#1\endcsname}
\def\label#1{\gdef\currentlabel{#1}}

\def\Lefteqlabel#1{\global\advance\eqnumer by 1
\label{#1}
\ifx\currentlabel\relax
\else
\expandafter\xdef
\csname\currentlabel\endcsname{(\the\numerosection.\the\eqnumer)}
\fi
\global\let\currentlabel\relax
\ifshowingMacros
\leqno\llap{%
\margefont #1\hphantom{M}}(\the\numerosection.\the\eqnumer)%
\else
\leqno(\the\numerosection.\the\eqnumer)%
\fi
}

\def\Righteqlabel#1{\global\advance\eqnumer by 1
\label{#1}
\ifx\currentlabel\relax
\else
\expandafter\xdef
\csname\currentlabel\endcsname{(\the\numerosection.\the\eqnumer)}%
\fi
\global\let\currentlabel\relax
\ifshowingMacros
\eqno(\the\numerosection.\the\eqnumer)%
\rlap{\margefont\hphantom{M}#1}
\else
\eqno(\the\numerosection.\the\eqnumer)%
\fi
}
\let\reqlabel=\Righteqlabel

\def\numlabel#1{\global\advance\eqnumer by1     
\label{#1}                                     
\ifx\currentlabel\relax                        
\else
\expandafter\xdef
\csname\currentlabel\endcsname{(\the\numerosection.\the\eqnumer)}
\fi
\global\let\currentlabel\relax
\ifshowingMacros
   (\the\numerosection.\the\eqnumer)\rlap{\margefont\hphantom{M}#1}
\else
(\the\numerosection.\the\eqnumer)%
\fi
}

\def\numero{\global\advance\eqnumer by1       
\number\numerosection.\number\eqnumer}

\def\lemlabel#1{\global\advance\lemnumer by 1
\ifshowingMacros%
   \marginalsigne{\margefont #1}%
\else
    \relax%\null
\fi
\label{#1}%
\ifx\currentlabel\relax%
\else
\expandafter\xdef%
\csname\currentlabel\endcsname{\the\numerosection.\the\lemnumer}%
\fi
\global\let\currentlabel\relax%
\the\numerosection.\the\lemnumer}

\def\numlem{\global\advance\lemnumer by1   
\the\numerosection.\number\lemnumer}  

%  ---------------------------------
%  Macros for automatic numbering of bibliography
%  ---------------------------------

\let\itemBibli=\item
\def\bibl#1#2\endbibl{\par{\itemBibli{#1} #2\par}}
\def\ref.#1.{{\csname#1\endcsname}}   
\newcount\bib   \bib=0
\def\bibmac#1{\advance\bib by 1       
\expandafter
\xdef\csname #1\endcsname{\the\bib}}
\def\BibMac#1{\advance\bib by1
\bibl{[\the\bib]}%
{\csname#1\endcsname}\endbibl}
\def\MakeBibliography#1{
\noindent{\bf #1}
\bigskip
\bib=0
\let\bibmac=\BibMac
\BiblioFil
\BiblioOrd
}

\newif\ifbouquin

\bouquinfalse
\font\tenssi=cmssi10 at 10 true pt
 at 10 true pt
\def\bibliostyle#1#2#3{%
{\rm #1}. \ifbouquin{\tenssi #2\/}\global\bouquinfalse%
                       \else{\it #2\/}\fi, {\rm #3}.}

%  --------------------------
%  The bibliograpy set 
%  ---------------------------

\def\BiblioFil{
\def\Calderon{\bibliostyle{\letrautor A.P.~Calder\'on}{On an inverse boundary value
problem}{Seminars on Numerical Analysis and Application to Continuum Physics,
SBM (Rio de Janeiro 1980), pp. 65--73}}
\def\Cessenatum{\bibliostyle{\letrautor M.~Cessenat}{Th\'eor\`emes de trace $L^p$ pour des espaces
de fonctions de la neutronique}{C.~R.~Acad.\ Sci.\ Paris, S\'erie I, {Vol.~299}
(1984), pp. 831--834}}
\def\Cessenatdois{\bibliostyle{\letrautor M.~Cessenat}{Th\'eor\`emes de trace pour des espaces
de fonctions de la neutronique}{C.~R.~Acad.\ Sci.\ Paris, S\'erie I, {Vol.~300}
(1985), pp. 89--92}}
\def\ChoulliStefanov{\bibliostyle{\letrautor M.~Choulli \& P.~Stefanov}{Inverse
scattering and inverse boundary value problem for the linear Boltzmann
equation}{Comm.\ Part.\ Diff.\ Equations, {Vol.~21} (5\&6), (1996), pp.~763--785}}
\def\CipMotRob{\bibliostyle{\letrautor R.~Cipolatti, C.M.~Motta \& N.C.~Roberty}{Stability Estimates for an 
Inverse Problem for the Linear Boltzmann Equation}{Revista Matem\'atica Complutense, Vol.~19, No.~1, 
(2006), pp. 113--132}}
\def\DautrayLions{\bibliostyle{\letrautor R.~Dautray \& J.-L.~Lions}{Mathematical
Analysis and Numerical Methods for Science and Technology}{Vol.~6, Springer-Verlag, 1993}}
\def\Folland{\bibliostyle{\letrautor G.B.~Folland}{Introduction to Partial Differential
Equations}{Mathematical No\-tes, 17, Princeton University Press, 1976}}
\def\Isakovdois{\bibliostyle{\letrautor V.~Isakov}{Inverse Problems for Partial Differential
Equations}{Applied Math.\ Sciences, Vol.~{127}, Springer, New York, 1998}}
\def\Rolci{\bibliostyle{\letrautor R.~Cipolatti}{Identification of the collision kernel in the
linear Boltzmann equation by a finite number of measurements on the boundary}{Comp.~Appl.~Math., 
Vol.~25, No.~2-3, (2006), pp.~331--351}}
\def\RolciIvo{\bibliostyle{\letrautor R.~Cipolatti \& Ivo F.~Lopez}{Determination of
coefficients for a dissipative wave equation via boundary measurements}{J.~Math.~Anal.\ 
Appl., Vol.~306, (2005), pp.~317--329}}
\def\Kharroubi{\bibliostyle{\letrautor M.~Mokhtar-Kharroubi}{Mathematical Topics in Neutron
Transport Theory - New Aspects}{Series on Advances in Mathematics
for Applied Sciences, Vol.~46, World Scientific, 1997}}
\def\ReedSimon{\bibliostyle{\letrautor M.~Reed \& B.~Simon}{Methods of Modern
Physics}{Vol.~3, Springer-Verlag, 1993}}
\def\Romanovum{\bibliostyle{\letrautor V.G.~Romanov}{Estimation of stability in the problem of 
determining the attenuation coefficient and the scattering indicatrix for the transport
equation}{Sibirsk.\ Mat.\ Zh.\ Vol.~37, No.~2, (1996), pp.~361--377, iii; translation in
Siberian Math.~J., Vol.~37, No.~2, (1996), pp.~308--324}}
\def\Romanovdois{\bibliostyle{\letrautor V.G.~Romanov}{Stability estimate in the three-domensional 
inverse problem for the transport equation}{J.~Inverse Ill-Posed Probl., Vol.~5, No.~5,
(1997), pp.~463--475}}
\def\StefUhlm{\bibliostyle{\letrautor P.~Stefanov \& G.~Uhlmann}{Optical tomography in two 
dimensions}{Methods and Applications of Analysis, Vol.~10, (2002), pp.~1--9}}
\def\Tamasan{\bibliostyle{\letrautor A.~Tamasan}{An inverse boundary value problem in two
dimensional transport}{Inverse Problems, {Vol.~18}, (2002), pp.~209--219}}
\def\Wang{\bibliostyle{\letrautor J.N.~Wang}{Stability estimates of an inverse problem for the
stationary transport equation}{Ann.~Inst.\ Henri Poincar\'e, {Vol.~70}, N.~5, (1999),
pp.~473--495}}
}

\def\BiblioOrd{
\bibmac{Calderon}
\bibmac{Cessenatum}
\bibmac{Cessenatdois}
\bibmac{ChoulliStefanov}
\bibmac{Rolci}
\bibmac{RolciIvo}
\bibmac{CipMotRob}
\bibmac{DautrayLions}
\bibmac{Folland}
\bibmac{Kharroubi}
\bibmac{ReedSimon}
\bibmac{Romanovum}
\bibmac{Romanovdois}
\bibmac{StefUhlm}
\bibmac{Tamasan}
\bibmac{Wang}
}
\BiblioOrd

%%%  ******************************************************
%%%                    BEGIN DOCUMENT 
%%%  ******************************************************

\centerline{\Titulo Identification of the Coefficients in the Linear Boltzmann Equation}
\smallskip  
\centerline{\Titulo by a Finite Number of Boundary Measurements}

\vskip 1truecm
\centerline{\letrautor Rolci Cipolatti}

\bigskip
{\baselineskip=14truept
\letranota
\centerline{Instituto de Matem\'atica}
\centerline{Universidade Federal do Rio de Janeiro}
\centerline{C.P.~68530, CEP 21945-970, Rio de Janeiro, Brazil}
\centerline{E-mail: cipolatti@im.ufrj.br}
}

\vskip 1truecm
{\parindent=1truecm \narrower\letranota\baselineskip=12truept
\noindent{\bf Abstract.} In this paper we consider an inverse problem for the time dependent 
linear Boltzmann equation. It concerns the identification of the coefficients  
via a finite number of measurements on the boundary. We prove that the total extinction coefficient 
and the collision kernel can be uniquely determined 
by at most $k$ measurements on the boundary, provided that these coefficients belong to a 
finite $k$-dimensional vector space.\par
\par
\noindent {\letrasl Key words:} inverse problem, linear Boltzmann equation, albedo operator,
            boundary measurements.\par

\par}
\vskip 1truecm
\noindent{\titulosecao \numsection. Introduction}
\bigskip
\noindent In this paper we consider an inverse problem 
 for the linear Boltzmann equation 
$$\partial_t u+\omega\cdot\nabla_xu+q(u-K_\kappa[u])=0
   \hbox{\rm\ \ in\ \ }(0,T)\times\S\times\Omega,\reqlabel{Boltzmann}$$
where $T>0$, $\Omega$ is a smooth bounded convex domain of $\R^N$, $N\ge2$, 
$\S$ denotes the unit sphere of $\R^N$, $q\in L^\infty(\Omega)$ 
and $K_\kappa$ is the integral operator with kernel 
$\kappa(x,\omega',\omega)$ defined by
$$K_\kappa[u](t,\omega,x)=\int_{\S}\kappa(x,\omega',\omega)u(t,\omega',x)\,d\omega'.$$

In applications, the equation \cite{Boltzmann} describes 
the dynamics of a monokinetic flow of particles in a body $\Omega$ under 
the assumption that the interaction between them is negligible. 
For instance, in the case of a low-density 
flux of neutrons (see [\cite{DautrayLions}], [\cite{ReedSimon}]), $q\ge 0$ is the total 
extinction coefficient (the inverse mean free path) and the collision kernel $\kappa$ is given by
$$\kappa(x,\omega',\omega)=c(x)h(x,\omega'\cdot\omega),$$
where $c$ corresponds to the within-group scattering probability and 
$h$ describes the anisotropy of the scattering process. In this model, 
$q(x)u(t,\omega,x)$ describes the loss of particles at $x$ in the direction $\omega$ at time $t$
due to absorption or scattering and $q(x)K_\kappa[u](t,\omega,x)$
represents the production of particles at $x$ 
in the direction $\omega$ from those coming from directions $\omega'$.  

Our focus here is the inverse problem of recovery the coefficients in \cite{Boltzmann}
via boundary measurements. More precisely, we are interested to recover $q$ and $\kappa$
by giving the incoming flux of particles on the boundary and measuring the outgoing one. 
Since these operations are described mathematically by the 
{\sl albedo\/} operator ${\cal A}_{q,\kappa}$,
%$${\cal A}_{q,\kappa}:L^1\bigl(0,T;L^1(\Sigma^-;d\xi)\bigr)\longrightarrow
%   L^1\bigl(0,T;L^1(\Sigma^+;d\xi)\bigr)$$ 
(which will be defined in the sequel), a first general mathematical question concerning
this inverse problem is to know if the knowledge of ${\cal A}_{q,\kappa}$ uniquely determines
$q,\kappa$, i.e., if the map $(q,\kappa)\mapsto {\cal A}_{q,\kappa}$ is invertible.  

We may precise this question. A first one is to know if the knowledge for ${\cal A}_{q,\kappa}(f)$
for all $f$ determines $(q,\kappa)$ ({\sl infinitely many measurementes\/});
a second one is to know if the knowledge of ${\cal A}_{q,\kappa}(f_j)$, 
for $j=1,2,\ldots,k$, determines $(q,\kappa)$ ({\sl finite number of measurements\/}).

There is a wide bibliography devoted to the first problem. We mention the
general results obtained by Choulli and Stefanov [\cite{ChoulliStefanov}]: 
they show that $q$ and $\kappa$ are uniquely 
determined by the {\sl albedo\/} operator (see also [\cite{Kharroubi}]). 
We also mention the stability results obtained 
by Cipolatti, Motta and Roberty (see [\cite{CipMotRob}] and the references therein).

There is also a lot of papers concerning the stationary case (see for instance those by
V.G.~Romanov [\cite{Romanovum}], [\cite{Romanovdois}], P.~Stefanov and G.~Uhlmann [\cite{StefUhlm}],
Tamasan [\cite{Tamasan}], J.N.~Wang [\cite{Wang}], and also the references therein).

In this work we focus on the second question, concerning the identification by a 
finite number of measurements. Under certain hypothesis and assuming that 
$\kappa(t,\omega',\omega)=c(x)h(\omega',\omega)$,
we prove that $q$ and $c$ can be uniquely determined by at most $k$ measurements, provided that $q$ and $c$
belongs to a finite $k$-dimensional vector subspace of $C(\overline\Omega)$.  To be more precise, we consider
the initial-boundary value problem
$$\left\{\eqalign{
&\partial_t u(t,\omega,x)+\omega\cdot\nabla u(t,\omega,x)+q(x)u(t,\omega,x) =
    q(x)K_\kappa[u](t,\omega,x), \cr
& u(0,\omega,x) = 0,\quad (\omega,x)\in \S\times\Omega, \cr
& u(t,\omega,\sigma) = f(t,\omega,\sigma),\quad (\omega,\sigma)
  \in \Sigma^-,\,\,\, t\in(0,T), \cr }\right.\reqlabel{IniBounProb}$$
where $f$ is the incoming flux and $\Sigma^-\colon=\{(\omega,\sigma)\in\S\times\partial\Omega
\,;\,\omega\cdot\nu(\sigma)<0\}$ is the incoming
part of the boundary. Then our main result can be stated as follows:
\smallskip
\noindent{\bf Theorem\ \lemlabel{MainTheorem}:} {\sl Let $\Omega\subset \R^N$ be a bounded convex domain of
class $C^1$, $M>0$, $T>\diam(\Omega)$ and ${\cal X}\colon=\hbox{\rm span}\{\rho_1,\rho_2,\ldots,\rho_k\}$, where 
$\{\rho_1,\rho_2,\ldots,\rho_k\}$ is a linearly independent subset of $C(\overline\Omega)$. 
We assume that $q\in L^\infty(\Omega)$, $\norma{q}{\infty}\le M$ and 
$\kappa\in L^\infty\bigl(\Omega;C(\S\times\S)\bigr)$. 

{\parindent=0.7truecm
\item{a)} If $q\in{\cal X}$, then there exist 
$\widetilde\omega_1,\ldots,\widetilde\omega_k\in\S$ and 
$f_1,\ldots,f_k\in C_0\bigl((0,T)\times\Sigma^-\bigr)$ 
that determine $q$ uniquely.\par 
\item{b)} If $\kappa(x,\omega',\omega)=c(x)h(\omega',\omega)$, 
where $c\in {\cal X}$ and $h\in C(\S\times\S)$ satisfies $h(\omega,\omega)\not=0$ for every $\omega\in\S$, 
then there exist $\widetilde\omega_1,\ldots,\widetilde\omega_k\in\S$ and 
$f_1,\ldots,f_k\in C_0\bigl((0,T)\times\Sigma^-\bigr)$
that determine $c$ uniquely.
\par}}

\smallskip          
\noindent{\bf Remark:} As we can see from the proof of Theorem~\cite{MainTheorem}, 
the functions $f_j$, $j=1,\ldots,k$, have the form
$$f_j(t,\omega,\sigma)\colon=\cases{
\phi_j(\sigma-t\omega)\ee^{\unim\lambda(t-\omega\cdot\sigma)} & in the case ({\sl a}),\cr
\delta_{\tilde\omega_j}(\sigma)\phi_j(\sigma-t\omega)\ee^{\unim\lambda(t-\omega\cdot\sigma)}& in the case ({\sl b}),\cr
}$$
where  $\lambda>0$, $t\in(0,T)$, $(\omega,\sigma)\in\Sigma^-$, $\phi_j\in C_0^\infty(\R^N\setminus\Omega)$ and 
$\delta_{\tilde\omega_j}$ is the spherical atomic measure 
concentrated on $\widetilde\omega_j$. The coefficients are identified by  
measuring the corresponding solutions on the outgoing part of the boundary, 
only in the directions $\widetilde\omega_1,\ldots,\widetilde\omega_k$.
 
The proof of Theorem~\cite{MainTheorem} is  based on the construction of
highly oscillatory solutions 
introduced in [\cite{CipMotRob}] and some arguments already used by the author in [\cite{RolciIvo}].

We organize the paper as follows: in Section~2 we recall the standard functional framework in which
the problem \cite{IniBounProb} is well posed and the 
albedo operator is defined; in Section~3, we introduce the highly oscillatory functions that will be
used, in Section~4, to prove Theorem~\cite{MainTheorem}.
     
\bigskip\goodbreak
\noindent{\titulosecao \numsection. The Functional Framework}
\bigskip

\noindent In this section we introduce the notation and we recall some well known  
results on the Transport Operator and the semigroup it generates (see [\cite{CipMotRob}] 
and the references therein for the proofs).

Let $\Omega\subset\R^N$ ($N\ge2$) be a convex and 
bounded domain of class $C^1$ and 
$\S$ the unit sphere of $\R^N$. We denote by 
$Q\colon=\S\times\Omega$ and $\Sigma$ its
boundary, i.e., $\Sigma\colon=\S\times{\partial\Omega}$.
For $p\in[1,+\infty)$ we consider the space $\Lpp$ 
with the usual norm 
$$\norma{u}{\Hpp}\colon=\left(\int_Q\mod{u(\omega,x)}^p\,dx\,d\omega\right)^{1/p},$$
where $d\omega$ denotes the surface measure on $\S$ 
associated to the Lebesgue measure in $\R^{N-1}$.

For each $u\in \Hpp$ we define $A_0u$ by
$$(A_0u)(\omega,x)\colon=\omega\cdot\nabla_x u(\omega,x)=\sum_{k=1}^N\omega_k
     {\partial u\over\partial x_k}(\omega,x),\quad\omega=(\omega_1,\ldots,\omega_N)$$
where the derivatives are taken in the sense of distributions in $\Omega$.

One checks easily that setting $\Wpp{p}\colon=\{u\in \Hpp\,;\, A_0u\in \Hpp\}$,
the operator $\bigl(A_0,\Wpp{p}\bigr)$ is a closed densely defined operator and
$\Wpp{p}$ with the graph norm is a Banach space. 

For every $\sigma\in{\partial\Omega}$, we denote $\nu(\sigma)$ the unit outward normal 
at $\sigma\in {\partial\Omega}$ and we consider the sets (respectively, the 
incoming and outgoing boundaries)
$$\Sigma^\pm \colon=\{(\omega,\sigma)\in \S\times{\partial\Omega}\,;\, 
    {}\pm\omega\cdot\nu(\sigma)>0\}.$$

In order to well define the albedo operator as a trace operator on the outgoing 
boundary, we consider $L^p(\Sigma^\pm;d\xi)$, where 
$d\xi\colon=\mod{\omega\cdot\nu(\sigma)}d\sigma d\omega$, and we introduce 
the spaces
$$\Wpptil{p}^\pm\colon=\bigl\{u\in\Wpp{p}\,;\, u_{|_{\Sigma^\pm}}\in L^p(\Sigma^\pm;d\xi)\bigr\},$$
which are Banach spaces if equipped with the norms
$$\norma{u}{\Wpptil{p}^\pm}\colon=\left(\norma{u}{\Wpp{p}}^p+
  \int_{\Sigma^\pm}\mod{\omega\cdot\nu(\sigma)}\mod{u(\omega,\sigma)}^p
     \,d\sigma d\omega\right)^{1/p}.$$
The next two lemmas concerne the continuity and surjectivity of the trace operators 
(see [\cite{Cessenatum}], [\cite{Cessenatdois}] and [\cite{CipMotRob}]):
$$\gamma_\pm:\Wpptil{p}^\pm\rightarrow L^p(\Sigma^\mp;d\xi),\quad
  \gamma_\pm(u)\colon=u_{|_{\Sigma^\mp}}. \reqlabel{DefNewLine}$$
 
\noindent{\bf Lemma\ \lemlabel{TracoCont2}:} {\sl Let $1\le p<+\infty$. 
Then there exists $C>0$ (depending only on $p$) such that
$$\int_{\Sigma^\mp}\mod{\omega\cdot\nu(\sigma)}\mod{u(\omega,\sigma)}^p
    \,d\sigma d\omega \le C\norma{u}{\Wpptil{p}^\pm}^p,
       \quad\forall u\in\Wpptil{p}^\pm.\reqlabel{DesigLado1Lado2}$$
Moreover, if $p>1$ and $1/p+1/p'=1$,
we have the Gauss identity   
$$\int_{Q}\divergente_x(uv\omega)\,dxd\omega=
     \int_{\Sigma}\omega\cdot\nu(\sigma)
	      u(\omega,\sigma)v(\omega,\sigma)\,d\sigma d\omega,\reqlabel{Gauss}$$
for all $u\in\Wpptil{p}^\pm$ and $v \in\Wpptil{p'}^\pm$.		  
}

As an immediate consequence of Lemma~\cite{TracoCont2}, we can introduce the space
$$\Wpptil{p}\colon=\bigl\{f\in\Wpp{p}\,;\, \int_\Sigma\mod{\omega\cdot\nu(\sigma)}
\mod{f(\omega,\sigma)}^p\,d\omega d\sigma<+\infty\bigr\}$$
an we have that $\Wpptil{p}^+=\Wpptil{p}^-=\Wpptil{p}$ with equivalent norms.  

\smallskip
\noindent{\bf Lemma\ \lemlabel{TracoSobre2}:} {\sl The trace operators $\gamma_\pm$ 
are surjective from $\Wpptil{p}^\pm$ onto $L^p(\Sigma^\mp;d\xi)$. 
More precisely, for each
$f\in L^p(\Sigma^\mp;d\xi)$, there exists $h\in\Wpptil{p}^\pm$ such that
$\gamma_{\pm}(h)=f$ and 
$$\norma{h}{\Wpptil{p}^\pm}\le C\norma{f}{L^p(\Sigma^\mp,d\xi)},$$
where $C>0$ is independent of $f$.}

We consider the operator $A: D(A)\rightarrow \Hpp$, defined by
$(Au)(\omega,x)\colon=\omega\cdot\nabla u(\omega,x)$,	 
with $D(A)\colon=\{u\in\Wpptil{p}\,;\, u|_{\Sigma^-}=0\}$. 

%For the reader's convenience we gather below a few more or less well known 
%results concerning the operator $(A,D(A))$ and the semigroup it generates. 

\smallskip\goodbreak
\noindent{\bf Theorem\ \lemlabel{m-acretivo}:} {\sl The operator $A$ 
is $m$-accretive in $\Hpp$, for $p\in[1,+\infty)$.}

\smallskip
\noindent{\bf Corollary\ \lemlabel{PrincMaxim}:} {\sl Let $f\in L^p(Q)$, $p\in[1,+\infty)$ and
assume that $u\in D(A)$ is a solution of $u+Au=f$. If $f\ge 0$ a.e.\ in $Q$, then
$u\ge 0$ a.e.\ in $Q$. In particular, it follows that}
$$\norma{u}{L^1(Q)}\le \norma{f}{L^1(Q)}.$$

It follows from Theorem~\cite{m-acretivo} and Corollary~\cite{PrincMaxim} that 
the operator $A$ generates a positive semigroup $\{U_0(t)\}_{t\ge0}$ 
of contractions acting on $\Hpp$. 

Let $q\in L^\infty(\Omega)$ and $\kappa: \Omega\times\S\times\S\rightarrow\R$ 
be a real measurable function satisfying
$$\left\{\eqalign{
\int_\S\mod{\kappa(x,\omega',\omega)}\,d\omega' & \le M_1  \hbox{\rm\ a.e.\ } 
  \Omega\times\S, \cr
\int_\S\mod{\kappa(x,\omega',\omega)}\,d\omega  & \le M_2  \hbox{\rm\ a.e.\ } 
  \Omega\times\S. \cr
}\right.\reqlabel{HipExtra}$$
Associated to these functions, we define the following continuous operators:
{\parindent=.7cm 
\item{1)} $B\in{\cal L}(\Hpp,\Hpp)$ defined by $B[u](\omega,x)\colon=q(x)u(\omega,x)$,
\item{2)} $K_\kappa[u](\omega,x)\colon=\int_\S\kappa(x,\omega',\omega)u(\omega',x)\,d\omega'$.
\par}
\noindent It follows from \cite{HipExtra} that $K_\kappa\in{\cal L}(\Hpp,\Hpp)$ $\forall p\in[1,+\infty)$ 
and (see [\cite{DautrayLions}]) 
$$\norma{K_\kappa[u]}{\Hpp}\le M_1^{1/p'}M_2^{1/p}\norma{u}{\Hpp}\le 
  \max\{M_1,M_2\}\norma{u}{\Hpp}.\reqlabel{PropLions}$$

The operator $A+B-K_\kappa: D(A)\rightarrow \Hpp$ generates a $c_0$-semigroup $\{U(t)\}_{t\ge 0}$
on $\Hpp$ satisfying 
$$\norma{U(t)}{{\cal L}}\le\hbox{\rm e}^{Ct},\quad C\colon=\norma{q^-}{\infty}+M_2.$$

We consider the initial-boundary 
value problem for the linear Boltzmann equation
$$\left\{\eqalign{
&\partial_t u(t,\omega,x)+\omega\cdot\nabla u(t,\omega,x)+q(x)u(t,\omega,x) = 
     qK_\kappa[u](t,\omega,x) \cr
& u(0,\omega,x) = 0,\quad (\omega,x)\in \S\times\Omega \cr
& u(t,\omega,\sigma) = f(t,\omega,\sigma),\quad (\omega,\sigma)
  \in \Sigma^-,\,\,\, t\in(0,T), \cr }\right.\reqlabel{ProbPar2}$$
where $q\in L^\infty(\Omega)$ and 
$$K_\kappa[u](t,\omega,x)=\int_\S\kappa(x,\omega',\omega)u(t,\omega',x)\,d\omega',$$
with $\kappa$ satisfying \cite{HipExtra}.

By the previous results, it follows that, for 
$f\in L^p\bigl(0,T;L^p(\Sigma^-,d\xi)\bigr)$, $p\in[1,+\infty)$,
there exists a unique solution  
$u\in C\bigl([0,T];\Wpptil{p}\bigr)\cap C^1\bigl([0,T];\Lpp\bigr)$
of \cite{ProbPar2}. This solution $u$ allows us to define the {\sl albedo\/} 
operator
$$\displaylines{
{\cal A}_{q,\kappa}:L^p\bigl(0,T; L^p(\Sigma^-,d\xi)\bigr)\rightarrow 
   L^p\bigl(0,T; L^p(\Sigma^+,d\xi)\bigr) \cr
{\cal A}_{q,\kappa}[f](t,\omega,\sigma)\colon=u(t,\omega,\sigma),\quad 
   (\omega,\sigma)\in\Sigma^+.\cr
}$$
\indent As a consequence of Lemmas~\cite{TracoCont2} and \cite{TracoSobre2},
${\cal A}_{q,\kappa}$ is a linear and bounded operator. 
%So, we denote by $\norma{{\cal A}_{q}}{p}$
%its norm. In order to simplify our notation, we will denote from now on
%$\norma{\ \ }{L^\pm_p}=\norma{\ \ }{L^p(0,T;L^p(\Sigma^\pm,d\xi))}$.

We also consider the following backward-boundary 
value problem, called the {\sl adjoint problem\/} of \cite{ProbPar2}:
$$\left\{\eqalign{
&\partial_t v(t,\omega,x)+\omega\cdot\nabla v(t,\omega,x)-q(x)v(t,\omega,x) 
      = -qK_\kappa^*[v](t,\omega,x) \cr
& v(T,\omega,x) = 0,\quad (\omega,x)\in \S\times\Omega \cr
& v(t,\omega,\sigma) = g(t,\omega,\sigma),\quad (\omega,\sigma)\in \Sigma^+, 
   t\in (0,T), \cr
}\right.\reqlabel{ProbParDual2}$$
where $g\in L^{p'}\bigl(0,T;L^{p'}(\Sigma^+,d\xi)\bigr)$, $p'\in[1,+\infty)$ and, 
$$K_\kappa^*[v](t,\omega',x)\colon=\int_\S\kappa(x,\omega',\omega)v(t,\omega,x)\,d\omega,$$ 
with the corresponding albedo operator ${\cal A}^*_{q,\kappa}$
$$\displaylines{
{\cal A}^*_{q,\kappa}: L^{p'}\bigl(0,T; L^{p'}(\Sigma^+,d\xi)\bigr)
  \rightarrow L^{p'}\bigl(0,T; L^{p'}(\Sigma^-,d\xi)\bigr) \cr
{\cal A}^*_{q,\kappa}[g](t,\omega,\sigma)\colon=v(t,\omega,\sigma),\quad 
   (\omega,\sigma)\in\Sigma^-.\cr
}$$ 
The operators ${\cal A}_{q,\kappa}$ and ${\cal A}^*_{q,\kappa}$ satisfy the following 
property:

\noindent{\bf Lemma\ \lemlabel{Dualidade2}:} {\sl Let 
$f\in L^p\bigl(0,T\,;\,L^p(\Sigma^-;d\xi)\bigr)$ and 
$g\in L^{p'}\bigl(0,T\,;\, L^{p'}(\Sigma^+;d\xi)\bigr)$, where $p,p'\in(1,+\infty)$ are such 
that $1/p+1/p'=1$. Then, we have}
$$ \displaylines{\qquad\qquad
\int_0^T\!\!\!\int_{\Sigma^-}\!(\omega\cdot\nu(\sigma))f(t,\omega,\sigma)
 {\cal A}^*_{q,\kappa}[g](t,\omega,\sigma)\,d\sigma d\omega dt= \hfill \cr
   \hfill= - \int_0^T\!\!\!\int_{\Sigma^+}\!(\omega\cdot\nu(\sigma))
     g(t,\omega,\sigma){\cal A}_{q,\kappa}[f](t,\omega,\sigma)\,d\sigma d\omega dt.
      \qquad\qquad\cr}$$

\noindent{\bf Proof:} It is a direct consequence of Lemma~\cite{TracoCont2}. 
Let $u(t,\omega,x)$ the solution of \cite{ProbPar2} with boundary 
condition $f$ and $v(t,\omega,x)$ the solution of \cite{ProbParDual2} with 
boundary $g$. We obtain the result by using \cite{Gauss}, once the equation in  
\cite{ProbPar2} is multiplied by $v$ and integrated over $Q$.\hfill\cqd

As a direct consequence of Lemma~\cite{Dualidade2}, we have:

\noindent{\bf Lemma\ \lemlabel{ParaDensid2}:} {\sl Let $T>0$, $q_{1},q_{2}\in 
L^\infty(\Omega)$ and 
$\kappa_1,\kappa_2$ satisfying \cite{HipExtra}. 
Assume that $u_1$ is the solution of \cite{ProbPar2} 
with coefficients $q_{1},\kappa_1$ and satisfying the boundary condition 
$f\in L^p\bigl(0,T;L^p(\Sigma^-,d\xi)\bigr)$, $p\in(1,+\infty)$ 
and that $u_2^*$ is the solution of \cite{ProbParDual2}, 
with $q_{2},\kappa_2$ and boundary condition
$g\in L^{p'}\bigl(0,T;L^{p'}(\Sigma^+,d\xi)\bigr)$, $1/p+1/p'=1$. 
Then we have }
$$\displaylines{
   \qquad\quad\int_0^T\!\!\!\int_{Q}\bigl(q_{2}(x)-q_{1}(x)\bigr)
     u_1(t,\omega,x)u_2^*(t,\omega,x)\,dxd\omega dt \hfill \cr
	\qquad\qquad\quad{}+\int_0^T\!\!\!\int_{Q}
       \bigl(q_{1}(x)K_{\kappa_1}[u_1](t,\omega,x)-q_2(x)K_{\kappa_2}[u_1](t,\omega,x)\bigr)
	       u_2^*(t,\omega,x)\,dxd\omega dt \hfill\cr
   \hfill{}=\int_0^T\!\!\!\int_{\Sigma^+}(\omega\cdot\nu(\sigma))
     \bigl[{\cal A}_{q_1,\kappa_1}[f]-{\cal A}_{q_2,\kappa_2}[f]\bigr]
       (t,\omega,\sigma)g(t,\omega,\sigma)\,d\sigma d\omega dt.\qquad\quad\cr	   
}$$

\goodbreak
\noindent{\titulosecao \numsection. Highly Oscillatory Solutions}
\bigskip
\noindent 
In this section we present some technical results related to special 
solutions of \cite{ProbPar2} and \cite{ProbParDual2}
that will be useful in the proof of Theorem~\cite{MainTheorem}.
They were published in [\cite{CipMotRob}] and [\cite{Rolci}], but for the reader's convenience
we present here their proofs. 
We denote by $\tilde q$ the zero extension of $q$ in the exterior of $\Omega$, i.e., 
$\widetilde q(x)=0$ for all $x\notin\Omega$.

\smallskip
\noindent{\bf Proposition\ \lemlabel{SolEsp2}:} {\sl Let $T>0$, 
$q_{1},q_{2}\in L^\infty(\Omega)$, and 
$\kappa$ satisfying \cite{HipExtra}.
Let $\psi_1,\psi_2\in C\bigl(\S,C^\infty_0(\R^N)\bigr)$ such that
$$\supp\psi_1(\omega,\cdot)\cap\overline\Omega =
  \left(\supp\psi_2(\omega,\cdot)+T\omega\right)\cap\overline\Omega =\emptyset,
      \quad\forall\omega\in\S.\reqlabel{CondSupp2}$$
Then, there exists $C_0>0$ such that, for each $\lambda>0$,  
there exist $R_{1,\lambda}\in C\bigl([0,T];\Wpptil{2}\bigr)$ and 
$R_{2,\lambda}^*\in C\bigl([0,T];\Wpptil{2}\bigr)$ satisfying
$$\norma{R_{1,\lambda}}{C([0,T];L^2(Q))}\le C_0,\quad
   \norma{R_{2,\lambda}^*}{C([0,T];L^{2}(Q))} \le C_0,\reqlabel{CondLim2}$$
for which the functions $u_1,u_2^*$ defined by
$$\left\{\eqalign{
u_1(t,\omega,x) & \colon= \psi_1(\omega,x-t\omega)\ee^{-\int_0^t\tilde q_{1}(x-s\omega)\,ds}
               \ee^{\unim\lambda(t-\omega\cdot x)}+R_{1,\lambda}(t,\omega,x)  \cr
u_2^*(t,\omega,x) & \colon= 
   \psi_2(\omega,x-t\omega)\ee^{\int_0^t\tilde q_{2}(x-s\omega)\,ds}
               \ee^{-\unim\lambda(t-\omega\cdot x)}+R_{2,\lambda}^*(t,\omega,x)  \cr
}\right.\reqlabel{SolEspeciais}$$
are solutions of \cite{ProbPar2} with $q=q_{1}$ and \cite{ProbParDual2} 
with $q=q_{2}$ respectively.
Moreover, if $\kappa\in L^\infty\bigl(\Omega;L^2(\S\times\S)\bigr)$, then we have
$$\lim_{\lambda\rightarrow+\infty}\norma{R_{1,\lambda}}{C([0,T];L^2(Q))}=
  \lim_{\lambda\rightarrow+\infty}\norma{R_{2,\lambda}^*}{C([0,T]L^2(Q))}=0.
\reqlabel{CondLim3}$$
}

\noindent{\bf Proof:} Let $u$ be the function  
$$u(t,\omega,x)\colon=\psi_1(\omega,x-t\omega)\ee^{-\int_0^t\tilde q_{1}(x-s\omega)\,ds}\,
         \ee^{\unim\lambda(t-\omega\cdot x)}+R(t,\omega,x).\reqlabel{Defu2}$$		  
By direct calculations, we easily verify that 
$$
\partial_tu+\omega\cdot\nabla u +q_{1}u -q_1K_{\kappa}[u] = 
     \partial_tR+\omega\cdot\nabla R +q_{1}R -q_1K_{\kappa}[R]-
        \ee^{\unim\lambda t}q_1Z_{1,\lambda},  
$$
where 
$$Z_{1,\lambda}(t,\omega,x)\colon=\int_\S\kappa(x,\omega',\omega)\psi_1(\omega',x-t\omega')
   \ee^{-\int_0^t\tilde q_{1}(x-s\omega')\,ds}\ee^{-\unim\lambda\omega'\cdot x}d\omega'.
\reqlabel{DefZ1}$$
By choosing $R_{1,\lambda}\in C^1\bigl([0,T];L^2(Q)\bigr)\cap C\bigl([0,T];D(A)\bigr)$ 
the solution of 
$$\left\{\eqalign{
\partial_tR+\omega\cdot\nabla R+q_{1}R  & = q_1K_{\kappa}[R]+\ee^{\unim\lambda t}q_1
   Z_{1,\lambda}, \cr
R(0,\omega,x)     & =0,  \quad (\omega,x)\in S\times\Omega,  \cr
R(t,\omega,\sigma) & =0,  \quad (\omega,\sigma)\in\Sigma^-,  \cr
}\right.\reqlabel{EqR2}$$ 
we see that \cite{CondSupp2} implies that the function $u$ defined 
by \cite{Defu2} satisfies \cite{ProbPar2} with boundary condition
$$f_\lambda^1(t,\omega,\sigma)\colon=
    \psi_1(\omega,\sigma-t\omega)\ee^{-\int_0^t\tilde q_{1}(\sigma-s\omega)\,ds}
          \ee^{\unim\lambda(t-\omega\cdot\sigma)},\quad(\omega,\sigma)\in\Sigma^-.$$
		  
Multiplying both sides of the equation in \cite{EqR2} by the complex conjugate 
of $R$, integrating it over $Q$ and taking its real part, we get, from 
Lemma~\cite{TracoCont2}, 
$$\displaylines{
\qquad{1\over 2}{d\hfil\over dt}\int_{Q}\mod{R(t)}^2d\omega dx+
   {1\over 2}\int_{\Sigma^+}\omega\cdot\nu(\sigma)\mod{R(t)}^2d\omega d\sigma+
     \int_Qq_{1}\mod{R(t)}^2d\omega dx-{}  \hfill\cr
\hfill\Re\int_Qq_1K_{\kappa}[R](t)\overline{R(t)}d\omega dx=
	   \Re\left[\ee^{\unim\lambda t}\int_Qq_1Z_{1,\lambda}(t)\overline{R(t)}\,
	     d\omega dx\right].\qquad\cr  
}$$  
It follows from the Cauchy-Schwarz inequality and \cite{PropLions} that
$$\int_Q\mod{K_{\kappa}[R(t)]}\mod{R(t)}\,dxd\omega\le 
   C_1\norma{R(t)}{L^2(Q)}^2,$$
where $C_1\colon=\max\{M_1,M_2\}$. Therefore, we obtain
$${d\hfil\over dt}\norma{R(t)}{L^2(Q)}^2\le
   C_2\norma{q_1}{\infty}\norma{R(t)}{L^2(Q)}^2+\norma{q_1}{\infty}\norma{Z_{1,\lambda}(t)}{L^2(Q)}^2, $$
where $C_2\colon=3+2C_1$. Since $R(0)=0$, we get, by integrating this last inequality on $[0,t]$, 
$$\eqalign{
 \norma{R(t)}{L^2(Q)}^2 & \le \norma{q_1}{\infty}\ee^{\norma{q_1}{\infty}TC_2}
     \int_0^t\norma{Z_{1,\lambda}(\tau)}{L^2(Q)}^{2}\,d\tau \cr
 & \le \norma{q_1}{\infty}\ee^{\norma{q_1}{\infty}TC_2}
   \norma{Z_{1,\lambda}}{L^2((0,T)\times Q)}^2\cr	 
 },\quad\forall t\in[0,T].\reqlabel{***1}
$$
	
The first inequality in \cite{CondLim2} follows easily because 
$\mod{Z_{1,\lambda}(t,\omega,x)}\le\norma{\psi_1}{\infty}\ee^{\norma{q_{1}}{\infty}T}M_1$
and, as the same arguments hold for $u_2^*$ and $R_{2,\lambda}^*$, 
we also obtain the second inequality.

We assume now $\kappa\in L^\infty\bigl(\Omega;L^2(\S\times\S)\bigr)$.
For each $x\in \R^N$, the map $\omega'\mapsto \exp(\unim\lambda\omega'\cdot x)$
converges weakly to zero in $L^2(\S)$ when $\lambda\rightarrow+\infty$ 
and the integral operator with kernel $\kappa(x,\cdot,\cdot)$ is compact 
in $L^2(\S)$. So, we obtain from \cite{DefZ1},
$$\lim_{\lambda\rightarrow+\infty}\norma{Z_{1,\lambda}(t,\cdot,x)}{L^2(\S)}=0
\quad\hbox{\rm a.e.\ in}\quad [0,T]\times\Omega.$$
Moreover, $\norma{Z_{1,\lambda}(t,\cdot,x)}{L^2(\S)}\le C$, where $C>0$ is a constant that 
does not depend on $\lambda$.
The Lebesgue's Dominated Convergence Theorem implies that
$$\lim_{\lambda\rightarrow+\infty}\norma{Z_{1,\lambda}}{L^2([0,T]\times Q)}=0.\reqlabel{***2}$$
From \cite{***2} and \cite{***1} we obtain \cite{CondLim3}, and our proof is complete. \hfill\cqd

\noindent{\bf Corollary\ \lemlabel{CorSolEsp2}:} {\sl Under the hypothesis of 
Proposition~\cite{SolEsp2}, if $q_1,q_2\in C(\overline{\Omega})$ and 
$\kappa\in L^\infty\bigl(\Omega;C(\S\times \S)\bigr)$, we have, for every $\omega\in \S$,
$$\lim_{\lambda\rightarrow+\infty}
    \norma{R_{1,\lambda}(\cdot,\omega,\cdot)}{C([0,T];L^2(\Omega))}=
  \lim_{\lambda\rightarrow+\infty}
    \norma{R_{2,\lambda}^*(\cdot,\omega,\cdot)}{C([0,T];L^2(\Omega))}=0.$$ 
}

\noindent{\bf Proof:} By multiplying both sides of the equation 
in \cite{EqR2} by the complex conjugate of $R(t,\omega,x)$,
integrating it over $\Omega$, taking its real part and applying the H\"older inequality, we get 
$$\eqalign{
 {d\hfil\over dt}\norma{R(t,\omega)}{L^2(\Omega)}^2 & \le
   4\norma{q_1}{\infty}\norma{R(t,\omega)}{L^2(\Omega)}^2\cr
   &\qquad{}+ \norma{q_1}{\infty}\left(\norma{K_\kappa[R](t,\omega)}{L^2(\Omega)}^2
  	+\norma{Z_{1,\lambda}(t,\omega)}{L^2(\Omega)}^2\right).\cr
}\reqlabel{***3}$$

Since 
$$\eqalign{
 \mod{K_\kappa[R](t,\omega,x)} & \le 
   \int_\S\mod{\kappa(x,\omega',\omega)}\mod{R(t,\omega',x)}\,d\omega'\cr
 &\le \left(\int_\S\mod{\kappa(x,\omega',\omega)}\,d\omega'\right)^{1/2}
      \left(\int_\S\mod{\kappa(x,\omega',\omega)}\mod{R(t,\omega',x)}^2\,d\omega'\right)^{1/2}\cr
 &\le M_1^{1/2}\norma{\kappa}{\infty}^{1/2}\left(\int_\S\mod{R(t,\omega',x)}^2\,d\omega'\right)^{1/2},\cr	   
}$$
we obtain
$$\norma{K_\kappa[R](t,\omega)}{L^2(\Omega)}^2\le M_1\norma{\kappa}{\infty}
   \norma{R(t)}{L^2(Q)}^2.\reqlabel{***4}$$

From \cite{***1}, \cite{***3} and \cite{***4} we have
$$\eqalign{
 {d\hfil\over dt}\norma{R(t,\omega)}{L^2(\Omega)}^2 & \le
   4\norma{q_1}{\infty}\norma{R(t,\omega)}{L^2(\Omega)}^2\cr
   &\qquad{}+ C\left(\norma{Z_{1,\lambda}}{L^2((0,T)\times Q)}^2
  	+\norma{Z_{1,\lambda}(t,\omega)}{L^2(\Omega)}^2\right).\cr
}$$   

Now, integrating this last inequality on time, we get
$$\eqalign{
 \norma{R(t,\omega)}{L^2(\Omega)}^2 & \le C\ee^{\norma{q_1}{\infty}T}
 \left(t\norma{Z_{1,\lambda}}{L^2((0,T)\times Q)}^2
  +\int_0^t\norma{Z_{1,\lambda}(\tau,\omega)}{L^2(\Omega)}^2\,d\tau\right)\cr
 &\le C\ee^{\norma{q_1}{\infty}T}\left(T\norma{Z_{1,\lambda}}{L^2((0,T)\times Q)}^2+
      \norma{Z_{1,\lambda}(\cdot,\omega,\cdot)}{L^2((0,T)\times\Omega)}^2\right).\cr
}$$

From Proposition~\cite{SolEsp2} we know that 
$\norma{Z_{1,\lambda}}{L^2((0,T)\times Q)}\rightarrow 0$ as $\lambda\rightarrow+\infty$.
On the other hand, as the map $\omega'\mapsto \ee^{\unim\omega'\cdot x}$ converges 
weakly to zero in $L^2(\S)$, we have from \cite{DefZ1}, for almost $x\in\Omega$,
$$\lim_{\lambda\rightarrow\infty}Z_{1,\lambda}(t,\omega,x)=0,\quad 
  \forall \omega\in S,\,\,\forall t\in[0,T]$$
and the conclusion follows from the Lebesgue's Theorem.\hfill\cqd

\noindent{\bf Lemma\ \lemlabel{Auxiliar}:} {\sl We assume that $q\in L^\infty(\Omega)$
and $\kappa$ satisfies \cite{HipExtra}. Let $S_{\lambda}^*$ the solution
of 
$$\left\{\eqalign{
\partial_tS+\omega\cdot\nabla S-qS  & = -qK_{\kappa}^*[S]+q\ee^{-\unim\lambda t}Z,\cr
S(T,\omega,x)     & =0,  \quad (\omega,x)\in \S\times\Omega,  \cr
S(t,\omega,\sigma) & =0,  \quad (\omega,\sigma)\in\Sigma^+,  \cr
}\right.\reqlabel{EqAux}$$
where $Z\in H^1\bigl(0,T;L^2(Q)\bigr)$ such that $Z(T)=0$. Then we have 
$$ \norma{S_{\lambda}^*}{C([0,T];L^2(Q))} \le C_0\quad\hbox{\sl and}\quad 
  \lim_{\lambda\rightarrow\infty}\norma{S_{\lambda}^*}{H^{-1}(0,T;L^2(Q))}=0,\reqlabel{ConcLemAux}$$ 
where $C_0$ is a constant independent of $\lambda$.}

\noindent{\bf Proof:} Multiplying both sides of the equation in \cite{EqAux} by the complex conjugate 
of $S_\lambda^*$, integrating it over $Q$ and taking its real part, we get 
$$\displaylines{
\quad{1\over 2}{d\hfil\over dt}\norma{S_\lambda^*(t)}{L^2(Q)}^2  
  +{1\over 2}\int_{\Sigma^-}(\omega\cdot\nu(\sigma))\mod{S_\lambda^*(t,\omega,\sigma)}\,d\omega d\sigma
   \ge {}-\norma{q}{\infty}\norma{S_\lambda^*(t)}{L^2(Q)}^2\hfill\cr
\hfill  {}-\norma{q}{\infty}\norma{K_\kappa^*[S](t)}{L^2(Q)}\norma{S_\lambda^*(t)}{L^2(Q)}-
  \norma{q}{\infty}\norma{Z(t)}{L^2(Q)}\norma{S_\lambda^*(t)}{L^2(Q)}\cr
}$$
Since $\norma{K_\kappa^*[S](t)}{L^2(Q)}\le \max\{M_1,M_2\}\norma{S_\lambda^*(t)}{L^2(Q)}$, we have
$${d\hfil\over dt}\norma{S_\lambda^*(t)}{L^2(Q)}^2\ge -C_2\norma{S_\lambda^*(t)}{L^2(Q)}^2-
    \norma{q}{\infty}\norma{Z(t)}{L^2(Q)}^2,$$
where $C_2\colon=(3+2\max\{M_1,M_2\})\norma{q}{\infty}.$
Integrating this last inequality on $[t,T]$ and taking into account that $S_\lambda^*(T)=0$, we obtain
$$
 \norma{S_\lambda^*(t)}{L^2(Q)}^2  \le  
  \norma{q}{\infty}\ee^{C_2T}\int_t^T\norma{Z(\tau)}{L^2(Q)}^2\,d\tau
   \le \norma{q}{\infty}\ee^{C_2T}\norma{Z}{L^2(0,T;L^2(Q))}\reqlabel{Aux-1}$$
and the inequality in \cite{ConcLemAux} follows easily.
 
We consider now
$$w_\lambda(t,\omega,x)\colon=\int_t^T S_{\lambda}^*(\tau,\omega,x)\,d\tau,\quad
  h(t,\omega,x)\colon=\int_t^T \ee^{-\unim\lambda\tau }Z(\tau,\omega,x)\,d\tau.\reqlabel{MudVarAux}$$
Then, it is easy to check that $w_\lambda$ satisfies
$$\left\{\eqalign{
\partial_tw+\omega\cdot\nabla w-qw  & = -qK_{\kappa}^*[w]+qh,\cr
w(T,\omega,x)     & =0,  \quad (\omega,x)\in \S\times\Omega,  \cr
w(t,\omega,\sigma) & =0,  \quad (\omega,\sigma)\in\Sigma^+,  \cr
}\right.\reqlabel{EqAuxbis}$$

Multiplying both sides of the equation in \cite{EqAuxbis} by the complex conjugate 
of $w_\lambda$, integrating it over $Q$, taking its real part and applying the Cauchy-Schwarz inequality, 
we get as before,
$$
 \norma{w_\lambda(t)}{L^2(Q)}^2  \le  \norma{q}{\infty}\ee^{C_2T}
 \norma{h}{L^2(0,T;L^2(Q))}^2 \le \norma{q}{\infty}T^2\ee^{C_2T}
    \norma{Z}{L^2(0,T;L^2(Q))}^2.\reqlabel{Aux-2}$$

As $S_\lambda^*=-\partial_t w_\lambda$, it follows from 
\cite{Aux-1} and \cite{Aux-2} that the set 
$\{w_{\lambda}\}$ is bounded in $C^1\bigl([0,T];L^2(Q)\bigr)$ 
and, in particular, is relatively compact in $C\bigl([0,T];L^2(Q)\bigr)$. 

On the other hand,  by integrating by parts the second integral in \cite{MudVarAux}, 
it is easy to check that there exists $C>0$ (depending only on $T$) such that 
$$\norma{h}{L^2(0,T:L^2(Q))}\le {C\over \mod{\lambda}}
    \norma{Z}{H^{1}(0,T;L^2(Q))}.\reqlabel{FormAux1}$$
Hence, by \cite{Aux-2}, it follows that $\norma{w_{\lambda}}{C([0,T];L^2(Q))}\rightarrow 0$
as $\lambda\rightarrow\infty$. Since the partial derivative in $t$,
$\partial_t:C\bigl([0,T];L^2(Q)\bigr) \rightarrow H^{-1}\bigl(0,T;L^2(Q)\bigr)$,
is a continuous operator, there exists a constant $C_3>0$ such that
$$\norma{S_{\lambda}^*}{H^{-1}(0,T;L^2(Q))}=\norma{\partial_t w_{\lambda}}{H^{-1}(0,T;L^2(Q))}
\le C_2\norma{w_{\lambda}}{C(0,T;L^2(Q))}$$ 
and we have the conclusion.\hfill\cqd

\bigskip\goodbreak
\noindent{\bf \numsection. \ Identification by a Finite Number of Boundary Measurements}
\bigskip

\noindent In this section we assume that  $\{\rho_1,\rho_2,\ldots,\rho_k\}$
is a given linearly independent set of functions of $C(\overline{\Omega})$ 
and we denote ${\cal X}\colon=\hbox{\rm span}\{\rho_1,\rho_2,\ldots,\rho_k\}$. 
For each $\widetilde\omega\in \S$ 
we consider $P_{\tilde\omega}[\rho_i]$ the X-ray transform of $\rho_i$ in the direction
$\widetilde\omega$, i.e., 
$$P_{\tilde\omega}[\rho_i](x)\colon=\int_{-\infty}^\infty\rho_i(x+t\widetilde\omega)\,dt$$
and, for each $\varepsilon>0$, 
$\Omega_\varepsilon\colon=\left\{x\in\R^N\setminus\overline\Omega\,;\, \dist(x,\Omega)<{\varepsilon}\right\}$.

The following Lemma, which the proof is given in [\cite{RolciIvo}], will be essential for the proof
of Theorem \cite{MainTheorem}:

\smallskip
\noindent{\bf Lemma\ \lemlabel{IndepDirec}:} {\sl For all $\varepsilon>0$, there exist 
$\widetilde\omega_j\in \S$ and $\phi_j\in C_0^\infty(\Omega_\varepsilon)$, $j=1,\ldots,k$,
such that the matrix ${A}=(a_{ij})$, with
entries defined by
$$a_{ij}\colon=\int_{\R^N}P_{\tilde\omega_j}[\rho_i](x)\phi_j^2(x)\,dx,\reqlabel{Entries}$$
is invertible.}

\noindent{\bf Remark\ \lemlabel{Obs}:} It follows from Lemma~\cite{IndepDirec}
and the equivalence of norms in finite dimensional vector spaces that there exists $C>0$ such
that, for all $\rho\in{\cal X}$,
$$\sum_{j=1}^k\Mod{\int_{\R^N}P_{\tilde\omega_j}[\rho](x)\phi_j^2(x)\,dx}\ge C\norma{\rho}{\infty}.$$

In order to prove Theorem \cite{MainTheorem}, we define, for $0<r<1$, the function
$\chi_r: \S\times \S\rightarrow\R$ as 
$\chi_r(\widetilde\omega,\omega)\colon=P(r\widetilde\omega,\omega)$, where $P$ 
is the Poisson kernel for $B_1(0)$, i.e., 
$$P(x,y)\colon={1-\mod{x}^2\over\alpha_N\mod{x-y}^N}.$$
From the well known properties of $P$ (see [\cite{Folland}]), we have 
$$\eqalign{
& \int_\S\chi_r(\widetilde\omega,\omega)\,d\omega
  =1,\,\,\,\forall\, r\in(0,1),\,\,\forall\widetilde\omega\in \S, \cr
& \lim_{r\rightarrow 1}\int_\S\chi_r(\widetilde\omega,\omega)\psi(\omega)\,d\omega
  =\psi(\widetilde\omega),\cr
}\reqlabel{PropPoisson}$$
where the above limit is taken in the topology of $L^p(\S)$, $p\in[1,+\infty)$ and uniformly on $\S$
if $\psi\in C(\S)$. We are now in position to prove Theorem~\cite{MainTheorem}.

\smallskip
\noindent{\bf Proof of Theorem \cite{MainTheorem}:} Let $\varepsilon\colon=(T-\diam(\Omega))/2$.

\noindent {\sl Step a\/}: We assume that $\kappa_1=\kappa_2=\kappa$ and $q_1,q_2\in{\cal X}$.
We define $\psi_1(\omega,x)=\varphi(x)$ and, for $\widetilde\omega\in \S$, 
$\psi_2(\omega,x)=\chi_r(\widetilde\omega,\omega)\psi(x)$, where
$0<r<1$ and $\varphi,\psi\in C^\infty_0(\Omega_\varepsilon)$, with $\norma{\psi}{L^\infty}\le 1$. 
Then $\psi_1$ and $\psi_2$
satisfy the condition \cite{CondSupp2} and we may consider the solutions $u_1$ and $u_2^*$
defined by \cite{SolEspeciais}, i.e., 
$$\eqalign{
u_1(t,\omega,x)&\colon=\varphi(x-t\omega)
  \ee^{-\int_0^t\tilde q_1(x-\tau\omega)d\tau}\ee^{\unim\lambda(t-x\cdot\omega)}
  +R_{1,\lambda}(t,\omega,x),\cr
u_2^*(t,\omega,x)&\colon=\chi_r(\widetilde\omega,\omega)\psi(x-t\omega)
  \ee^{\int_0^t\tilde q_2(x-\tau\omega)d\tau}\ee^{-\unim\lambda(t-x\cdot\omega)}
  +R_{2,\lambda,r}^*(t,\omega,x),\cr  
}$$  
where  $\lambda>0$ will be chosen a posteriori. We shall write 
$$\eqalign{
 \Phi_{1,\lambda}(t,\omega,x)&\colon=
   \varphi(x-t\omega)\ee^{-\int_0^t\tilde q_1(x-\tau\omega)d\tau}\ee^{\unim\lambda(t-x\cdot\omega)} \cr
  \Psi_{2,\lambda}(t,\omega,x)&\colon=
   \psi(x-t\omega)\ee^{\int_0^t\tilde q_2(x-\tau\omega)d\tau}\ee^{-\unim\lambda(t-x\cdot\omega)}  \cr 
}$$
in such a way that 
$$\eqalign{
u_1(t,\omega,x) & =\Phi_{1,\lambda}(t,\omega,x)+R_{1,\lambda}(t,\omega,x),\cr
u_2^*(t,\omega,x) & =\chi_r(\widetilde\omega,\omega)
  \Psi_{2,\lambda}(t,\omega,x)+ R_{2,\lambda,r}^*(t,\omega,x).\cr
}$$
  
Substituting $u_1$ and $u_2^*$ in the identity given in Lemma~\cite{ParaDensid2}, we have
$$I(\lambda,r)-J(\lambda,r)=L(\lambda,r),\reqlabel{ParFazLim1-r}$$
where 
$$\eqalign{
I(\lambda,r)&\colon=\int_0^T\!\!\!\int_{Q}\bigl(q_{2}(x)-q_{1}(x)\bigr)
     u_1(t,\omega,x)u_2^*(t,\omega,x)\,dxd\omega dt, \cr 
J(\lambda,r)&\colon=\int_0^T\!\!\!\int_Q\bigl(q_2(x)-q_1(x)\bigr)K_\kappa[u_1](t,\omega,x)
   u_2^*(t,\omega,x)\,dxd\omega dt,\cr
L(\lambda,r)&\colon=\int_0^T\!\!\!\int_{\Sigma^+}(\omega\cdot\nu(\sigma))\bigl(
    {\cal A}_1[f_{\lambda}]-{\cal A}_2[f_{\lambda}]\bigr)g_{2,\lambda,r}\,d\sigma d\omega dt.\cr
}$$
In the above formulas, we are denoting ${\cal A}_i={\cal A}_{q_i}$, $i=1,2$ and  
$$\eqalign{	   
f_{\lambda}(t,\omega,\sigma)&\colon=\varphi(\sigma-t\omega)\ee^{-\int_0^t\tilde q_1(\sigma-\tau\omega)d\tau}
  \ee^{\unim\lambda(t-\sigma\cdot\omega)},\quad(\omega,\sigma)\in\Sigma^-,\cr
g_{2,\lambda,r}(t,\omega,\sigma)&\colon=\chi_r(\widetilde\omega,\omega)\psi(\sigma-t\omega)
  \ee^{\int_0^t\tilde q_2(\sigma-\tau\omega)d\tau}\ee^{-\unim\lambda(t-\sigma\cdot\omega)},
   \quad(\omega,\sigma)\in\Sigma^+.\cr   
}$$
Since $\Omega$ is convex, for $s>0$ and $(\omega,\sigma)\in\Sigma^-$, we have $(\sigma-s\omega)\notin\Omega$ and
$\widetilde q_1(\sigma-s\omega)=0$. Hence $f_\lambda$ does not depend on $q_1$, i.e.,
$$f_{\lambda}(t,\omega,\sigma)=\varphi(\sigma-t\omega)
  \ee^{\unim\lambda(t-\sigma\cdot\omega)},\quad(\omega,\sigma)\in\Sigma^-.\reqlabel{NDef-fg}$$

By denoting $\rho(x)=\widetilde q_2(x)-\widetilde q_1(x)$ 
and by considering the special form of 
$u_1$ and $u_2^*$, we may write $I(\lambda,r)$ and $J(\lambda,r)$
as $I=I_1+I_2+I_3+I_4$ and $J=J_1+J_2+J_3+J_4$, where
$$\eqalign{
I_1(\lambda,r)&\colon=\int_0^T\!\!\!\int_Q\rho(x)\ee^{-\int_0^t\rho(x-s\omega)ds}\chi_r(\widetilde\omega,\omega)
    \varphi(x-t\omega)\psi(x-t\omega)\,dxd\omega dt,\cr
I_2(\lambda,r)&\colon=\int_0^T\!\!\!\int_Q\rho(x)\Phi_{1,\lambda}(t,\omega,x)R_{2,\lambda,r}^*(t,\omega,x)\,dxd\omega dt,\cr
I_3(\lambda,r)&\colon=\int_0^T\!\!\!\int_Q\rho(x)\chi_r(\widetilde\omega,\omega)
    \Psi_{2,\lambda}(t,\omega,x)R_{1,\lambda}(t,\omega,x)\,dxd\omega dt,\cr   
I_4(\lambda,r)&\colon=\int_0^T\!\!\!\int_Q\rho(x)R_{1,\lambda}(t,\omega,x)R_{2,\lambda,r}^*(t,\omega,x)\,dxd\omega dt.\cr	
}$$
and
$$\eqalign{
J_1(\lambda,r)&\colon=\int_0^T\!\!\!\int_Q\rho(x)\left[\int_\S \kappa(x,\omega',\omega)
    \Phi_{1,\lambda}(t,\omega',x)\,d\omega'\right]\chi_r(\widetilde\omega,\omega)\Psi_{2,\lambda}(t,\omega,x)\,dxd\omega dt,\cr
J_2(\lambda,r)&\colon=\int_0^T\!\!\!\int_Q\rho(x)\left[\int_\S \kappa(x,\omega',\omega)
    \Phi_{1,\lambda}(t,\omega',x)\,d\omega'\right]R_{2,\lambda,r}^*(t,\omega,x)\,dxd\omega dt,\cr
J_3(\lambda,r)&\colon=\int_0^T\!\!\!\int_Q\rho(x)\left[\int_\S \kappa(x,\omega',\omega)
    R_{1,\lambda}(t,\omega',x)\,d\omega'\right]\chi_r(\widetilde\omega,\omega)
	\Psi_{2,\lambda}(t,\omega,x)\,dxd\omega dt,\cr   
J_4(\lambda,r)&\colon=\int_0^T\!\!\!\int_Q\rho(x)\left[\int_\S \kappa(x,\omega',\omega)
    R_{1,\lambda}(t,\omega',x)\,d\omega'\right]R_{2,\lambda,r}^*(t,\omega,x)\,dxd\omega dt.\cr	
}$$

Taking the limit as $r\rightarrow 1^-$ in the above expressions, we get from \cite{PropPoisson},
$I_i(\lambda,r)\rightarrow I_i(\lambda)$ and $J_i(\lambda,r)\rightarrow J_i(\lambda)$, $i=1,\ldots,4$, where
$$\eqalign{
I_1(\lambda)&\colon=\int_0^T\!\!\!\int_\Omega\rho(x)\ee^{-\int_0^t\rho(x-s\tilde\omega)ds}
    \varphi(x-t\widetilde\omega)\psi(x-t\widetilde\omega)\,dxdt,\cr
I_2(\lambda)&\colon=\int_0^T\!\!\!\int_Q\rho(x)\Phi_{1,\lambda}(t,\omega,x)S^*_{2,\lambda}(t,\omega,x)\,dxd\omega dt,\cr	
I_3(\lambda)&\colon=\int_0^T\!\!\!\int_\Omega\rho(x)
    \Psi_{2,\lambda}(t,\widetilde\omega,x)R_{1,\lambda}(t,\widetilde\omega,x)\,dxdt,\cr   
I_4(\lambda)&\colon=\int_0^T\!\!\!\int_Q\rho(x)R_{1,\lambda}(t,\omega,x)S^*_{2,\lambda}(t,\omega,x)\,dxd\omega dt,\cr
}$$	

$$\eqalign{
J_1(\lambda)&\colon=\int_0^T\!\!\!\int_\Omega \rho(x)\left[\int_\S \kappa(x,\omega',\widetilde\omega)
    \Phi_{1,\lambda}(t,\omega',x)d\omega'\right]
	\Psi_{2,\lambda}(t,\widetilde\omega,x)\,dxdt,\cr
J_2(\lambda)&\colon=\int_0^T\!\!\!\int_Q \rho(x)\left[\int_\S \kappa(x,\omega',\omega)
	\Phi_{1,\lambda}(t,\omega',x)d\omega'\right]S^*_{2,\lambda}(t,\omega,x)\,dxd\omega dt,\cr
J_3(\lambda)&\colon=\int_0^T\!\!\!\int_\Omega \rho(x)\left[\int_\S \kappa(x,\omega',\widetilde\omega)
    R_{1,\lambda}(t,\omega',x)\,d\omega'\right]\Psi_{2,\lambda}(t,\widetilde\omega,x)\,dxdt,\cr 
J_4(\lambda)&\colon=\int_0^T\!\!\!\int_Q\rho(x)\left[\int_\S \kappa(x,\omega',\omega)
    R_{1,\lambda}(t,\omega',x)\,d\omega'\right]S^*_{2,\lambda}(t,\omega,x)\,dxd\omega dt,\cr	 
}$$
and $S_{2,\lambda}^*$ is the unique solution of 
$$\left\{\eqalign{
\partial_tS+\omega\cdot\nabla S-q_2(S-K_{\kappa}^*[S]) & = \ee^{-\unim\lambda t}q_2
   Z_{2,\lambda}^*, \cr
S(T,\omega,x)     & =0,  \quad (\omega,x)\in \S\times\Omega,  \cr
S(T,\omega,\sigma) & =0,  \quad (\omega,\sigma)\in\Sigma^+,  \cr
}\right.\reqlabel{EqR2-lim}$$
where 
$$Z_{2,\lambda}^*(t,\omega,x)\colon=\kappa(x,\omega,\widetilde\omega)\psi(x-t\widetilde\omega)
   \ee^{\int_0^t\tilde q_2(x-s\tilde\omega)ds}\ee^{\unim\lambda x\cdot\tilde\omega}.$$
Moreover, from \cite{NDef-fg} and \cite{PropPoisson}, it follows that 
$L(\lambda,r)\rightarrow L(\lambda)$, where
$$L(\lambda)\colon  = \int_0^T\!\!\!\int_{\partial\Omega}(\widetilde\omega\cdot\nu(\sigma))^+
    \bigl(\widetilde{\cal A}_1[f_{\lambda}]-\widetilde{\cal A}_2[f_{\lambda}]\bigr)(t,\widetilde\omega,\sigma)
	  \Psi_{2,\lambda}(t,\widetilde\omega,\sigma)\,d\sigma dt,\reqlabel{NMF1}$$
where $\widetilde{\cal A}_i[f_{\lambda}]$ denotes the zero extension of 
${\cal A}_i[f_{\lambda}]$ on $\partial\Omega$.
Therefore, by taking the limit as $r\rightarrow 1^-$ in \cite{ParFazLim1-r}, we have
$$\sum_{i=1}^4I_i(\lambda)-\sum_{i=1}^4J_i(\lambda)=L(\lambda).$$	  
So,
$$\mod{I_1(\lambda)}\le \sum_{i=2}^4\mod{I_i(\lambda)}+\sum_{i=1}^4\mod{J_i(\lambda)}+\mod{L(\lambda)}\reqlabel{ParaEstim}$$

In what follows, we prove that there exists $C(\lambda)>0$, $C(\lambda)\rightarrow 0$ as $\lambda\rightarrow\infty$, such that 
$$\sum_{i=2}^4\mod{I_i(\lambda)}+\sum_{i=1}^4\mod{J_i(\lambda)}\le C(\lambda)\norma{\rho}{\infty}.\reqlabel{Clambda}$$

We begin with terms of odd indices. Remembering that $\norma{\psi}{\infty}\le 1$, it is easy to see that
$$\eqalign{
 \mod{I_3(\lambda)} & \le C\norma{\rho}{\infty}
    \norma{R_{1,\lambda}(\cdot,\widetilde\omega,\cdot)}{C(0,T;L^2(\Omega))},\cr
 \mod{J_1(\lambda)} & \le C\norma{\rho}{\infty}
    \int_0^T\!\!\!\int_\Omega\mod{K_\kappa[\Phi_{1,\lambda}](t,\widetilde\omega,x)}\,dxdt\cr
 \mod{J_3(\lambda)} &\le C\norma{\rho}{\infty}
    \int_0^T\!\!\!\int_\Omega\mod{K_\kappa[R_{1,\lambda}](t,\widetilde\omega,x)}\,dxdt,\cr
}\reqlabel{desigualdade1}$$
where $C=C(T,M)$.
From Corollary~\cite{CorSolEsp2} it follows that 
$$\lim_{\lambda\rightarrow\infty}\norma{R_{1,\lambda}(\cdot,\widetilde\omega,\cdot)}{C(0,T;L^2(\Omega))}=0.\reqlabel{estim1}$$
Since $\omega'\mapsto \ee^{\unim\lambda x\cdot\omega'}$ converges weakly to zero in $L^2(\S)$ 
as $\lambda\rightarrow\infty$, we have $\Phi_{1,\lambda}(t,\cdot,x)\rightharpoonup 0$ a.e.\ $t$ and $x$.
Since $K_\kappa$ is a compact operator in $L^2(\S)$, we have
$K_\kappa[\Phi_{1,\lambda}](t,\widetilde\omega,x)\rightarrow 0$ a.e $t$ and $x$. So, by the Lebesgue Theorem, it follows that
$$\lim_{\lambda\rightarrow\infty}\int_0^T\!\!\!\int_\Omega\mod{K_\kappa[\Phi_{1,\lambda}](t,\widetilde\omega,x)}\,dxdt=0.\reqlabel{estim2}$$
Since
$$\mod{K_\kappa[R_{1,\lambda}](t,\widetilde\omega,x)}\le\int_\S\mod{\kappa(x,\omega',\widetilde\omega)R_{1,\lambda}(t,\omega',x)}\,d\omega'
\le\norma{\kappa}{\infty}\int_\S\mod{R_{1,\lambda}(t,\omega',x)}\,d\omega',$$
it follows that 
$$\mod{J_3(\lambda)}\le C\norma{\rho}{\infty}\norma{\kappa}{\infty}\norma{R_{1,\lambda}}{C(0,T;L^2(Q))}$$
and we have from Lemma~\cite{SolEsp2} 
$$\lim_{\lambda\rightarrow\infty}\norma{R_{1,\lambda}}{C(0,T;L^2(Q))}=0.\reqlabel{estim4}$$

On the other hand, we have for the even indices:
$$\eqalign{
\mod{I_2(\lambda)} &\le C\norma{\rho}{\infty}
    \int_0^T\!\!\!\int_Q\mod{\varphi(x-t\omega)S_{2,\lambda}^*(t,\omega,x)}\,dxd\omega dt,\cr
\mod{I_4(\lambda)} &\le C\norma{\rho}{\infty}\norma{R_{1,\lambda}}{L^2(0,T;L^2(Q))}
    \norma{S_{2,\lambda}^*}{L^2(0,T;L^2(Q))},\cr
\mod{J_2(\lambda)} &\le C\norma{\rho}{\infty}
    \int_0^T\!\!\!\int_Q\mod{K_\kappa[\Phi_{1,\lambda}](t,\omega,x)S_{2,\lambda}^*(t,\omega,x)}\,dxd\omega dt,\cr
\mod{J_4(\lambda)} &\le C\norma{\rho}{\infty}
   	\int_0^T\!\!\!\int_Q\mod{K_\kappa[R_{1,\lambda}](t,\omega,x)S_{2,\lambda}^*(t,\omega,x)}\,dxd\omega dt,\cr
}\reqlabel{desigualdade2}$$

From Proposition~\cite{SolEsp2} and Lemma~\cite{Auxiliar}, it follows that
$$\lim_{\lambda\rightarrow\infty}\norma{R_{1,\lambda}}{L^2(0,T;L^2(Q))}
    \norma{S_{2,\lambda}^*}{L^2(0,T;L^2(Q))}=0.\reqlabel{estim5}$$
Moreover, since $S_{2,\lambda}^*$ is bounded in $L^2(0,T;L^2(Q))$, $K_\kappa[\Phi_{1,\lambda}]\rightarrow 0$
and $K_\kappa[R_{1,\lambda}]\rightarrow 0$ in $L^2(0,T;L^2(Q))$,
and we have 
$$\eqalign{
 \lim_{\lambda\rightarrow\infty}\int_0^T\!\!\!\int_Q\mod{K_\kappa[\Phi_{1,\lambda}]
   (t,\omega,x)S_{2,\lambda}^*(t,\omega,x)}\,dxd\omega dt & =0, \cr
 \lim_{\lambda\rightarrow\infty}\int_0^T\!\!\!\int_Q\mod{K_\kappa[R_{1,\lambda}]
   (t,\omega,x)S_{2,\lambda}^*(t,\omega,x)}\,dxd\omega dt & =0.\cr
 }\reqlabel{estim6}$$ 

On the other hand, since $\varphi\in C_0^\infty(\Omega_\varepsilon)$, it follows from the choice of $\varepsilon$ that the map
$(t,\omega,x)\mapsto \varphi(x-t\omega)$ belongs to $H_0^1(0,T;L^2(Q))$. Therefore, 
by Lemma~\cite{Auxiliar},
$$\int_0^T\!\!\!\int_Q\mod{\varphi(x-t\omega)S_{2,\lambda}^*(t,\omega,x)}\,dxd\omega dt\le
\norma{\varphi}{H_0^1(0,T;L^2(Q))}\norma{S_{2,\lambda}^*}{H^{-1}(0,T;L^2(Q))}\rightarrow 0\reqlabel{estim6}$$
as $\lambda\rightarrow\infty$ and we conclude from \cite{desigualdade1}--\cite{estim6} that there exists $C(\lambda)$
satisfying \cite{Clambda}. So, it follows from \cite{ParaEstim} that
$$\mod{I_1(\lambda)}\le C(\lambda)\norma{\rho}{\infty}+
   \int_0^T\!\!\!\int_{\partial\Omega}\Mod{\bigl(\widetilde{\cal A}_1[f_\lambda]-
     \widetilde{\cal A}_2[f_\lambda]\bigr)(t,\widetilde\omega,\sigma)}\,dtd\sigma.\reqlabel{BoaEstim1}$$

Now we remark that
$$\eqalign{
\mod{I_1(\lambda)} & = \Mod{\int_0^T\!\!\!\int_\Omega\rho(x)
  \ee^{-\int_0^t\rho(x-s\tilde\omega)ds}\varphi(x-t\tilde\omega)\psi(x-t\tilde\omega)dxdt} \cr
                   & = \Mod{\int_{\R^N}\left[1-\ee^{-\int_0^T\rho(y+s\tilde\omega)ds}\right]\varphi(y)\psi(y)dy}.
}$$
By taking the supremum on $\psi\in L^\infty(\R^N)$, $\norma{\psi}{\infty}\le 1$, we get
$$\sup_{\psi}\Mod{\int_{\R^N}\left[1-\ee^{-\int_0^T\rho(y+s\tilde\omega)ds}\right]\varphi(y)\psi(y)dy}=
   \int_{\R^N}\Mod{1-\ee^{-\int_0^T\rho(y+s\tilde\omega)ds}}\mod{\varphi(y)}\,dy.$$
   
Since 
$$\Mod{1-\ee^{-\int_0^T\rho(y+s\tilde\omega)ds}}\ge \Mod{\int_0^T\rho(y+s\widetilde\omega)ds}\ee^{-MT},$$
we have from \cite{BoaEstim1} (with $C_0=\ee^{-MT}$)
$$\eqalign{
 C_0\Mod{\int_{\R^N}\!\int_0^T\rho(y+s\widetilde\omega)\mod{\varphi(y)}\,dsdy} & \le C(\lambda)\norma{\rho}{\infty}\cr
   &{}+\int_0^T\!\!\!\int_{\partial\Omega}\Mod{\bigl(\widetilde{\cal A}_1[f_\lambda]-
     \widetilde{\cal A}_2[f_\lambda]\bigr)(t,\widetilde\omega,\sigma)}\,dtd\sigma.\cr}
\reqlabel{BoaEstim2}$$   

Since $(\supp\varphi+s\widetilde\omega)\cap\Omega=\emptyset$ for all $|s|\ge T$, we have
$$\eqalign{
 \Mod{\int_{\R^N}\int_0^T\rho(y+s\widetilde\omega)\mod{\varphi(y)}\,dsdy} 
 &=\Mod{\int_{-\infty}^{\infty}\!\int_{\R^N}\rho(x)\mod{\varphi(x-s\tilde\omega)}\,dxds} \cr
 &=\Mod{\int_{\R^N}\int_{-\infty}^{\infty}\rho(y+s\widetilde\omega)\mod{\varphi(y)}\,dsdy} \cr
 &=\Mod{\int_{\R^N}P_{\tilde\omega}[\rho](y)\mod{\varphi(y)}\,dy}\cr
}\reqlabel{BoaEstim2}$$

From \cite{BoaEstim2} and \cite{BoaEstim1}, we get
$$\eqalign{
 C_0\Mod{\int_{\R^N}P_{\tilde\omega}[\rho](y)\mod{\varphi(y)}\,dy} & \le C(\lambda)\norma{\rho}{\infty}\cr
   &{}+\int_0^T\!\!\!\int_{\partial\Omega}\Mod{\bigl(\widetilde{\cal A}_1[f_\lambda]-
     \widetilde{\cal A}_2[f_\lambda]\bigr)(t,\widetilde\omega,\sigma)}\,dtd\sigma.\cr}
\reqlabel{BoaEstim3}$$

We are now in position to conclude the proof of step~a. First of all, we consider in \cite{BoaEstim3}
$\widetilde\omega$ as the directions $\widetilde\omega_1,\ldots,\widetilde\omega_k$ and $\varphi$ as the functions 
$\phi_1^2,\ldots,\phi_k^2$ given by Lemma~\cite{IndepDirec},
in such a way that (see \cite{NDef-fg})
$$f_{j,\lambda}(t,\omega,\sigma)=\phi_j^2(\sigma-t\omega)\ee^{\unim\lambda(t-\sigma\cdot\omega)}.$$
Hence, by Lemma~\cite{IndepDirec} and Remark~\cite{Obs}, we can write, for some $C>0$,
$$\eqalign{
 C\norma{q_1-q_2}{\infty} & \le C(\lambda)\norma{q_1-q_2}{\infty}\cr
   &{}+\sum_{j=1}^k\int_0^T\!\!\!\int_{\partial\Omega}\Mod{\bigl(\widetilde{\cal A}_1[f_{j,\lambda}]-
     \widetilde{\cal A}_2[f_{j,\lambda}]\bigr)(t,\widetilde\omega_j,\sigma)}\,dtd\sigma.\cr
}$$
Therefore, if ${\cal A}_1[f_{j,\lambda}](t,\widetilde\omega_j,\sigma)={\cal A}_2[f_{j,\lambda}](t,\widetilde\omega_j,\sigma)$
on $\Sigma^+_{\tilde\omega_j}$, for $j=1,\ldots,k$, we have
$$C\norma{q_1-q_2}{\infty}  \le C(\lambda)\norma{q_1-q_2}{\infty}$$
and the conclusion follows easily if we choose $\lambda>0$ large enough.

\medskip

\noindent {\sl Step b\/}:
We assume that $q_1=q_2=q$ and $\kappa_i(x,\omega',\omega)=
c_i(x)h(\omega',\omega)$, where $c_1, c_2\in{\cal X}$. For $\widetilde\omega\in \S$, we define 
$\psi_1(\omega,x)=\chi_s(\widetilde\omega,\omega)\phi(x)$ and 
$\psi_2(\omega,x)=\chi_r(\widetilde\omega,\omega)\phi(x)$, where
$0<r,s<1$ and $\phi\in C^\infty_0(\Omega_\varepsilon)$. 
Then $\psi_1$ and $\psi_2$
satisfy the condition \cite{CondSupp2} and we may consider the solutions $u_1$ and $u_2^*$
defined by \cite{SolEspeciais}, i.e., 
$$\eqalign{
u_1(t,\omega,x)&\colon=\chi_s(\widetilde\omega,\omega)\phi(x-t\omega)
  \ee^{-\int_0^t\tilde q(x-\tau\omega)d\tau}\ee^{\unim\lambda(t-x\cdot\omega)}
  +R_{1,\lambda,s}(t,\omega,x),\cr
u_2^*(t,\omega,x)&\colon=\chi_r(\widetilde\omega,\omega)\phi(x-t\omega)
  \ee^{\int_0^t\tilde q(x-\tau\omega)d\tau}\ee^{-\unim\lambda(t-x\cdot\omega)}
  +R_{2,\lambda,r}^*(t,\omega,x),\cr  
}$$  
where  $\lambda>0$ will be chosen a posteriori. We shall write 
$$\eqalign{
 \Phi_{\lambda}(t,\omega,x)&\colon=
   \phi(x-t\omega)\ee^{-\int_0^t\tilde q(x-\tau\omega)d\tau}\ee^{-\unim\lambda x\cdot\omega} \cr
  \Psi_{\lambda}(t,\omega,x)&\colon=
   \phi(x-t\omega)\ee^{\int_0^t\tilde q(x-\tau\omega)d\tau}\ee^{\unim\lambda x\cdot\omega}  \cr 
}$$
in such a way that 
$$\eqalign{
u_1(t,\omega,x) & =\ee^{\unim\lambda t}\chi_s(\widetilde\omega,\omega)
  \Phi_{\lambda}(t,\omega,x)+R_{1,\lambda,s}(t,\omega,x),\cr
u_2^*(t,\omega,x) & =\ee^{-\unim\lambda t}\chi_r(\widetilde\omega,\omega)
  \Psi_{\lambda}(t,\omega,x)+ R_{2,\lambda,r}^*(t,\omega,x).\cr
}\reqlabel{Defu1v2-rs}$$
  
Substituting $u_1$ and $u_2^*$ in the identity given in Lemma~\cite{ParaDensid2}, we have
$$J(\lambda,r,s)=L(\lambda,r,s),\reqlabel{ParFazLim-r}$$
where 
$$\eqalign{
J(\lambda,r,s)&\colon=\int_0^T\!\!\!\int_Qq(x)\bigl(c_1(x)-c_2(x)\bigr)K_h[u_1](t,\omega,x)
   u_2^*(t,\omega,x)\,dxd\omega dt,\cr
L(\lambda,r,s)&\colon=\int_0^T\!\!\!\int_{\Sigma^+}(\omega\cdot\nu(\sigma))\bigl(
    {\cal A}_1[f_{\lambda,s}]-{\cal A}_2[f_{\lambda,s}]\bigr)g_{\lambda,r}\,d\sigma d\omega dt.\cr
}$$
In the above formulas, we are denoting ${\cal A}_i={\cal A}_{c_i}$, $i=1,2$ and  
$$\eqalign{	   
f_{\lambda,s}(t,\omega,\sigma)&\colon=\chi_s(\widetilde\omega,\omega)\Phi_{\lambda}(t,\omega,\sigma),
  \quad(\omega,\sigma)\in\Sigma^-,\cr
g_{\lambda,r}(t,\omega,\sigma)&\colon=\chi_r(\widetilde\omega,\omega)\Psi_{\lambda}(t,\omega,\sigma),
   \quad(\omega,\sigma)\in\Sigma^+.\cr   
}\reqlabel{NDef-fg}$$
In particular, it follows from the definition of the Albedo Operator and \cite{Defu1v2-rs},
$${\cal A}_1[f_{\lambda,s}]-{\cal A}_2[f_{\lambda,s}]=
  R_{1,\lambda,s}-R_{2,\lambda,s}, \quad\hbox{\rm on\ }(0,T)\times \Sigma^+.\reqlabel{Dif-dos-Alb}$$ 
  
By denoting $\eta(x)=\widetilde q(x)\bigl(\widetilde c_1(x)-\widetilde c_2(x)\bigr)$ 
and by considering the special form of 
$u_1$ and $u_2^*$, we may write $J(\lambda,r,s)$
as $J=J_1+J_2+J_3+J_4,$ where
$$\eqalign{
J_1(\lambda,r,s)&\colon=\int_0^T\!\!\!\int_Q\eta(x)\left[\int_\S h(\omega',\omega)
    \chi_s(\widetilde\omega,\omega')\Phi_{\lambda}(t,\omega',x)d\omega'\right]\times{}\cr
	&\qquad\qquad{}\times	\chi_r(\widetilde\omega,\omega)\Psi_{\lambda}(t,\omega,x)\,dxd\omega dt,\cr
J_2(\lambda,r,s)&\colon=\int_0^T\!\!\!\int_Q\eta(x)\left[\int_\S h(\omega',\omega)
    \chi_s(\widetilde\omega,\omega')\Phi_{\lambda}(t,\omega',x)d\omega'\right]
	R_{2,\lambda,r}^*(t,\omega,x)\,dxd\omega dt,\cr
J_3(\lambda,r,s)&\colon=\int_0^T\!\!\!\int_Q\eta(x)\left[\int_\S h(\omega',\omega)
    R_{1,\lambda,s}(t,\omega',x)d\omega'\right]\chi_r(\widetilde\omega,\omega)
	\Psi_{\lambda}(t,\omega,x)\,dxd\omega dt,\cr   
J_4(\lambda,r,s)&\colon=\int_0^T\!\!\!\int_Q\eta(x)\left[\int_\S h(\omega',\omega)
    R_{1,\lambda,s}(t,\omega',x)d\omega'\right]R_{2,\lambda,r}^*(t,\omega,x)\,dxd\omega dt.\cr	
}$$

Taking the limit as $r\rightarrow 1^-$ in the above expressions, we get from \cite{IndepDirec}
$J_i(\lambda,r,s)\rightarrow J_i(\lambda,s)$, where 
$$\eqalign{
J_1(\lambda,s)&\colon=\int_0^T\!\!\!\int_\Omega \eta(x)\left[\int_\S h(\omega',\widetilde\omega)
    \chi_s(\widetilde\omega,\omega')\Phi_{\lambda}(t,\omega',x)d\omega'\right]
	\Psi_{\lambda}(t,\widetilde\omega,x)\,dxdt,\cr
J_2(\lambda,s)&\colon=\int_0^T\!\!\!\int_Q \eta(x)\left[\int_\S h(\omega',\omega)
    \chi_s(\widetilde\omega,\omega')\Phi_{\lambda}(t,\omega',x)d\omega'\right]
	S_{2,\lambda}^*(t,\omega,x)\,dxd\omega dt,\cr	
J_3(\lambda,s)&\colon=\int_0^T\!\!\!\int_\Omega \eta(x)\left[\int_\S h(\omega',\widetilde\omega)
    R_{1,\lambda,s}(t,\omega',x)d\omega'\right]\Psi_{\lambda}(t,\widetilde\omega,x)\,dxdt,\cr  
J_4(\lambda,s)&\colon=\int_0^T\!\!\!\int_Q\eta(x)\left[\int_\S h(\omega',\omega)
    R_{1,\lambda,s}(t,\omega',x)d\omega'\right]S_{2,\lambda}^*(t,\omega,x)\,dxd\omega dt\cr	 		 	
}$$
and $S_{2,\lambda}^*$ is the unique solution of 
$$\left\{\eqalign{
\partial_tS+\omega\cdot\nabla S-qS  & = -qK_{\kappa_2}^*[S]+\ee^{-\unim\lambda t}q
   Z_{2,\lambda}^*, \cr
S(T,\omega,x)     & =0,  \quad (\omega,x)\in \S\times\Omega,  \cr
S(T,\omega,\sigma) & =0,  \quad (\omega,\sigma)\in\Sigma^+,  \cr
}\right.\reqlabel{EqR2-lim}$$
  
Moreover, from \cite{NDef-fg} and \cite{PropPoisson}, it follows that 
$L(\lambda,r,s)\rightarrow L(\lambda,s)$, where
$$\eqalign{
 L(\lambda,s)\colon & =\int_0^T\!\!\!\int_{\partial\Omega}(\widetilde\omega\cdot\nu(\sigma))^+
    \bigl(\widetilde{\cal A}_1[f_{\lambda,s}]-\widetilde{\cal A}_2[f_{\lambda,s}]\bigr)(t,\widetilde\omega,\sigma)
	  \Psi_{\lambda}(t,\widetilde\omega,\sigma)\,d\sigma dt\cr
  &	=\int_0^T\!\!\!\int_{\partial\Omega}(\widetilde\omega\cdot\nu(\sigma))^+
    \bigl(R_{1,\lambda,s}(t,\widetilde\omega,\sigma)-R_{2,\lambda,s}(t,\widetilde\omega,\sigma)\bigr)
	  \Psi_{\lambda}(t,\widetilde\omega,\sigma)\,d\sigma dt,\cr 
}	   \reqlabel{NMF1}$$
where $\widetilde{\cal A}_i[f_{\lambda,i}]$ denotes the zero extension of 
${\cal A}_i[f_{\lambda,i}]$ on $\partial\Omega$.
Therefore, by taking the limit as $r\rightarrow 1^-$ in \cite{ParFazLim-r}, we have
$$J_1(\lambda,s)+J_2(\lambda,s)+J_3(\lambda,s)+J_4(\lambda,s)=L(\lambda,s).$$	  
	  
Now, it is time to take the limit as $s\rightarrow 1^-$. For the first two terms of the 
right hand side of the above identity, we get (for $i=1,2$)
$J_i(\lambda,s)\rightarrow J_i(\lambda)$, where
$$\eqalign{
J_1(\lambda)&\colon=\int_0^T\!\!\!\int_\Omega \eta(x)h(\widetilde\omega,\widetilde\omega)
    \Phi_{\lambda}(t,\widetilde\omega,x)\Psi_{\lambda}(t,\widetilde\omega,x)\,dxdt\cr
	&=h(\widetilde\omega,\widetilde\omega)\int_0^T\!\!\!\int_\Omega \eta(x)\phi(x-t\widetilde\omega)^2dxdt,\cr
J_2(\lambda)&\colon=\int_0^T\!\!\!\int_Q \eta(x)h(\widetilde\omega,\omega)
    \Phi_{\lambda}(t,\widetilde\omega,x)S_{2,\lambda}^*(t,\omega,x)\,dxd\omega dt.\cr
}\reqlabel{NM1}$$

On the other hand, the dependence on $s$ in the other terms is given by 
$R_{1,\lambda,s}$ and $R_{2,\lambda,s}$, which are the solution of ($j=1,2$)
$$\left\{\eqalign{
\partial_tR+\omega\cdot\nabla R +qR  & = qK_{\kappa_j}[R]+\ee^{\unim\lambda t}q
   Z_{j,\lambda,s}, \cr
R(0,\omega,x)     & =0,  \quad (\omega,x)\in \S\times\Omega,  \cr
R(0,\omega,\sigma) & =0,  \quad (\omega,\sigma)\in\Sigma^+,  \cr
}\right.\reqlabel{EqR2-r}$$  
where 
$$Z_{j,\lambda,s}(t,\omega,x)\colon=\int_\S\kappa_j(x,\omega',\omega)\chi_s(\widetilde\omega,\omega')
   \Phi_{\lambda}(t,\omega',x)d\omega'.\reqlabel{DefZ-2r}$$
It is an immediate consequence of \cite{PropPoisson} and the Lebesgue's Theorem that,  
as $s\rightarrow 1$, $Z_{j,\lambda,s}\rightarrow Z_{j,\lambda}$  in 
$C\bigl([0,T];L^2(Q)\bigr)$,  where
$$Z_{j,\lambda}(t,\omega,x)\colon=
   \kappa_j(x,\widetilde\omega,\omega)\Phi_{\lambda}(t,\widetilde\omega,x).\reqlabel{Def-Z2lim}$$
Hence, 
$$\lim_{s\rightarrow 1^-}R_{j,\lambda,s}= S_{j,\lambda}\quad 
\hbox{\rm in}\quad C\bigl([0,T];L^2(Q)\bigr),$$ 
where $S_{j,\lambda}$ is the solution of
$$\left\{\eqalign{
\partial_tS+\omega\cdot\nabla S+qS  & = qK_{\kappa_j}[S]+\ee^{\unim\lambda t}qZ_{j,\lambda}, \cr
S(0,\omega,x)     & =0,  \quad (\omega,x)\in \S\times\Omega,  \cr
S(t,\omega,\sigma) & =0,  \quad (\omega,\sigma)\in\Sigma^-,  \cr
}\right.\reqlabel{EqR1-limbis}$$
and $Z_{j,\lambda}(t,\omega,x)\colon=
   c_j(x)h(\widetilde\omega,\omega)\Phi_{\lambda}(t,\widetilde\omega,x)$. 
Therefore, $J_i(\lambda,s)\rightarrow J_i(\lambda)$, ($i=3,4$) and 
$L(\lambda,s)\rightarrow L(\lambda)$,  where
$$\eqalign{
J_3(\lambda)&\colon=\int_0^T\!\!\!\int_\Omega \eta(x)\left[\int_\S h(\omega',\widetilde\omega)
    S_{1,\lambda}(t,\omega',x)d\omega'\right]\Psi_{\lambda}(t,\widetilde\omega,x)\,dxdt,\cr
J_4(\lambda)&\colon=\int_0^T\!\!\!\int_Q \eta(x)\left[\int_\S h(\omega',\widetilde\omega)
    S_{1,\lambda}(t,\omega',x)d\omega'\right]S_{2,\lambda}^*(t,\omega,x)\,dxd\omega dt.\cr
L(\lambda)&\colon=\int_0^T\!\!\!\int_{\partial\Omega}\bigl(\widetilde\omega\cdot\nu(\sigma)\bigr)^+
\bigl(S_{1,\lambda}(t,\widetilde\omega,\sigma)-S_{2,\lambda}(t,\widetilde\omega,\sigma)\bigr)
    \Psi_\lambda(t,\widetilde\omega,\sigma)\,d\sigma dt	
}\reqlabel{NM2}$$
and we obtain
$$\mod{J_1(\lambda)}\le \mod{J_2(\lambda)}+\mod{J_3(\lambda)}+
   \mod{J_4(\lambda)}+\mod{L(\lambda)},\reqlabel{FimFim1}$$
where
$$\eqalign{
\mod{J_2(\lambda)} & \le \norma{\eta}{\infty}\norma{h}{\infty}\ee^{MT}
  \int_0^T\!\!\!\int_Q\mod{\phi(x-t\widetilde\omega)S_{2,\lambda}^*(t,\omega,x)}\,dxd\omega dt,\cr
\mod{J_3(\lambda)} & \le \norma{\eta}{\infty}\norma{\phi}{\infty}\ee^{MT}
  \norma{K_h[S_{1,\lambda}]}{L^2(0,T;L^2(Q))},\cr  
\mod{J_4(\lambda)} & \le \norma{\eta}{\infty}\norma{K_h[S_{1,\lambda}]}{L^2(0,T;L^2(Q))}
  \norma{S_{2,\lambda}^*}{L^2(0,T;L^2(Q))},\cr
\mod{L(\lambda)} & \le \norma{\phi}{\infty}\ee^{MT}\int_0^T\!\!\!\int_{\partial\Omega}
  \bigl(\widetilde\omega\cdot\nu(\sigma)\bigr)^+
  \Mod{S_{1,\lambda}(t,\widetilde\omega,\sigma)-S_{2,\lambda}(t,\widetilde\omega,\sigma)}\,d\sigma dt.       
}\reqlabel{FimFim2}$$

Since $\phi\in C_0^\infty(\Omega_\varepsilon)$, it follows from the choice of $\varepsilon$ that the function 
$(t,\omega,x)\mapsto\phi(x-t\widetilde\omega)$ 
belongs to $H_0^1(0,T;L^2(Q))$ (as a constant function on $\omega$). Hence, we have 
$$\mod{J_2(\lambda)}\le \norma{\rho}{\infty}\ee^{MT}
    \norma{\phi}{H_0^1(0,T;L^2(Q))}\norma{S_{2,\lambda}^*}{H^{-1}(0,T;L^2(Q))}.$$
On the other hand, from the weak convergence to zero in $L^2\bigl(0,T;L^2(Q)\bigr)$ of 
$S_{1,\lambda}$, it follows that 
$$\lim_{\lambda\rightarrow +\infty}\norma{K_h[S_{1,\lambda}]}{L^2(0,T;L^2(Q))}=0.\reqlabel{FimFim3}$$ 
Hence, we have from \cite{FimFim1}--\cite{FimFim3} and Lemma~\cite{Auxiliar}, 
$$\eqalign{
 \mod{J_1(\lambda)} & =\mod{h(\widetilde\omega,\widetilde\omega)}
  \Mod{\int_0^T\!\!\!\int_\Omega\eta(x)\phi(x-t\widetilde\omega)^2\,dxdt}\cr
  & \le C(\lambda)\norma{\eta}{\infty}+
     C_2\int_0^T\!\!\!\int_{\partial\Omega}\bigl(\widetilde\omega\cdot\nu(\sigma)\bigr)^+
	  \Mod{S_{1,\lambda}(t,\widetilde\omega,\sigma)-
       S_{2,\lambda}(t,\widetilde\omega,\sigma)}\,d\sigma dt,\cr
 }\reqlabel{FimFim4}$$
where $C(\lambda)\rightarrow 0$ as $\lambda\rightarrow+\infty$.

Since $(\supp\phi+s\widetilde\omega)\cap\Omega=\emptyset$ for all $|s|\ge T$, we have
$$\eqalign{
 \Mod{\int_0^T\!\!\!\int_\Omega\eta(x)\phi(x-t\widetilde\omega)^2\,dxdt} 
   & =\Mod{\int_{\R^N}\int_0^T\eta(y+s\widetilde\omega)\phi(y)^2\,dsdy} \cr
   & =\Mod{\int_{\R^N}\int_{-\infty}^{\infty}\rho(y+s\widetilde\omega)\phi(y)^2\,dsdy} \cr
   & =\Mod{\int_{\R^N}P_{\tilde\omega}[\eta](y)\phi(y)^2}\,dy\cr
}\reqlabel{FimFim5}$$
and we get
$$\displaylines{
\qquad \qquad\mod{h(\widetilde\omega,\widetilde\omega)}\Mod{\int_{\R^N}P_{\tilde\omega}[\eta](y)\phi(y)^2}\,dy
  \le C(\lambda)\norma{\eta}{\infty}+{}\hfill\cr
    \hfill C_2\int_0^T\!\!\!\int_{\partial\Omega}\bigl(\widetilde\omega\cdot\nu(\sigma)\bigr)^+
	 \Mod{S_{1,\lambda}(t,\widetilde\omega,\sigma)-
       S_{2,\lambda}(t,\widetilde\omega,\sigma)}\,d\sigma dt\qquad\qquad \cr
 } $$
 
We are now in position to conclude the proof. First of all, we consider in the above inequality 
the directions $\widetilde\omega_1,\ldots,\widetilde\omega_k$ and the functions 
$\phi_1,\ldots,\phi_k$ given by Lemma~\cite{IndepDirec},
in such a way that, from Remark~\cite{Obs}, we can write
$$\eqalign{
 C_0\norma{c_1-c_2}{\infty}&\le C(\lambda)\norma{c_1-c_2}{\infty}+{}\hfill\cr
   &{}+C_2\sum_{j=1}^k\int_0^T\!\!\!\int_{\partial\Omega}\bigl(\widetilde\omega\cdot\nu(\sigma)\bigr)^+
     \Mod{S_{1,\lambda}(t,\widetilde\omega_j,\sigma)-
        S_{2,\lambda}(t,\widetilde\omega_j,\sigma)}\,d\sigma dt,\cr
 }\reqlabel{ParaConcluir}$$
for some constant $C_0>0$. If we denote by
$$u_{ij}(t,\omega,\sigma)=\chi_s(\widetilde\omega_j,\omega)\Phi_\lambda(t,\omega,x)+
  R_{i,\lambda,s}(t,\omega,x),\quad i=1,2,\,\,j=1,\ldots,k,$$  
it follows from \cite{PropPoisson} that, as $s\rightarrow 1^-$, $u_{ij}\rightarrow u_{ij}^\#$, where 
$$u_{ij}^\#=\delta_{\tilde\omega_j}\Phi_\lambda + S_{i,\lambda},\quad i=1,2,\,\,j=1,\ldots,k$$ 
and $\delta_{\tilde\omega_j}$ is the spherical atomic measure concentrated on $\widetilde\omega_j$.
It is clear that $u_{1j}^\#-u_{2j}^\#= S_{1,\lambda}-S_{2,\lambda}$.  
Therefore, if $u_{1j}^\#(t,\widetilde\omega_j,\sigma)=u_{2j}^\#(t,\widetilde\omega_j,\sigma)$
on $\Sigma_{\tilde\omega_j}^+$, for $j=1,\ldots,k$, it follows from \cite{ParaConcluir} that 
$$C_0\norma{c_1-c_2}{\infty}\le C(\lambda)\norma{c_1-c_2}{\infty}$$
and the conclusion follows easily if we choose $\lambda>0$ small enough.\hfill\cqd

\bigskip\goodbreak
{\parindent=16pt
\MakeBibliography{\titulosecao REFERENCES}

}
\bye